\documentclass[]{article}
\usepackage{amssymb}
\font\black=cmbx10 \font\sblack=cmbx7 \font\ssblack=cmbx5
\font\blackital=cmmib10  \skewchar\blackital='177 \font\sblackital=cmmib7
\skewchar\sblackital='177 \font\ssblackital=cmmib5
\skewchar\ssblackital='177 \font\sanss=cmss10 \font\ssanss=cmss8 scaled
900 \font\sssanss=cmss8 scaled 600 \font\blackboard=msbm10
\font\sblackboard=msbm7 \font\ssblackboard=msbm5 \font\caligr=eusm10
\font\scaligr=eusm7 \font\sscaligr=eusm5 
\font\fraktur=eufm10 \font\sfraktur=eufm7 \font\ssfraktur=eufm5

\font\bsymb=cmsy10 scaled\magstep2
\def\all#1{\setbox0=\hbox{\lower1.5pt\hbox{\bsymb
       \char"38}}\setbox1=\hbox{$_{#1}$} \box0\lower2pt\box1\;}
\def\exi#1{\setbox0=\hbox{\lower1.5pt\hbox{\bsymb \char"39}}
       \setbox1=\hbox{$_{#1}$} \box0\lower2pt\box1\;}

\def\tx#1{{\fam0\relax#1}}

\newfam\bifam
\textfont\bifam=\blackital \scriptfont\bifam=\sblackital
\scriptscriptfont\bifam=\ssblackital
\def\bi#1{{\fam\bifam\relax#1}}

\newfam\blfam
\textfont\blfam=\black \scriptfont\blfam=\sblack
\scriptscriptfont\blfam=\ssblack

\newfam\bbfam
\textfont\bbfam=\blackboard \scriptfont\bbfam=\sblackboard
\scriptscriptfont\bbfam=\ssblackboard

\newfam\ssfam
\textfont\ssfam=\sanss \scriptfont\ssfam=\ssanss
\scriptscriptfont\ssfam=\sssanss
\def\sss#1{{\fam\ssfam\relax#1}}

\newfam\clfam
\textfont\clfam=\caligr \scriptfont\clfam=\scaligr
\scriptscriptfont\clfam=\sscaligr

\newfam\frfam
\textfont\frfam=\fraktur \scriptfont\frfam=\sfraktur
\scriptscriptfont\frfam=\ssfraktur

\def\hpb#1{\setbox0=\hbox{${#1}$}
    \copy0 \kern-\wd0 \kern.2pt \box0}
\def\vpb#1{\setbox0=\hbox{${#1}$}
    \copy0 \kern-\wd0 \raise.08pt \box0}

\def\pmb#1{\setbox0\hbox{${#1}$} \copy0 \kern-\wd0 \kern.2pt \box0}
\def\pmbb#1{\setbox0\hbox{${#1}$} \copy0 \kern-\wd0
      \kern.2pt \copy0 \kern-\wd0 \kern.2pt \box0}
\def\pmbbb#1{\setbox0\hbox{${#1}$} \copy0 \kern-\wd0
      \kern.2pt \copy0 \kern-\wd0 \kern.2pt
    \copy0 \kern-\wd0 \kern.2pt \box0}
\def\pmxb#1{\setbox0\hbox{${#1}$} \copy0 \kern-\wd0
      \kern.2pt \copy0 \kern-\wd0 \kern.2pt
      \copy0 \kern-\wd0 \kern.2pt \copy0 \kern-\wd0 \kern.2pt \box0}
\def\pmxbb#1{\setbox0\hbox{${#1}$} \copy0 \kern-\wd0 \kern.2pt
      \copy0 \kern-\wd0 \kern.2pt
      \copy0 \kern-\wd0 \kern.2pt \copy0 \kern-\wd0 \kern.2pt
      \copy0 \kern-\wd0 \kern.2pt \box0}

\mathchardef\za="710B  %\alpha
\mathchardef\zb="710C  %\beta
\mathchardef\zg="710D  %\gamma
\mathchardef\zd="710E  %\delta
\mathchardef\zve="710F %\epsilon
\mathchardef\zz="7110  %\zeta
\mathchardef\zh="7111  %\eta
\mathchardef\zvy="7112 %\theta
\mathchardef\zi="7113  %\iota
\mathchardef\zk="7114  %\kappa
\mathchardef\zl="7115  %\lambda
\mathchardef\zm="7116  %\mu
\mathchardef\zn="7117  %\nu
\mathchardef\zx="7118  %\xi
\mathchardef\zp="7119  %\pi
\mathchardef\zr="711A  %\rho
\mathchardef\zs="711B  %\sigma
\mathchardef\zt="711C  %\tau
\mathchardef\zu="711D  %\upsilon
\mathchardef\zvf="711E %\phi
\mathchardef\zq="711F  %\chi
\mathchardef\zc="7120  %\psi
\mathchardef\zw="7121  %\omega
\mathchardef\ze="7122  %\varepsilon
\mathchardef\zy="7123  %\vartheta
\mathchardef\zvp="7124 %\varpi
\mathchardef\zvr="7125 %\varrho
\mathchardef\zvs="7126 %\varsigma
\mathchardef\zf="7127  %\varphi
\mathchardef\zG="7000  %\Gamma
\mathchardef\zD="7001  %\Delta
\mathchardef\zY="7002  %\Theta
\mathchardef\zL="7003  %\Lambda
\mathchardef\zX="7004  %\Xi
\mathchardef\zP="7005  %\Pi
\mathchardef\zS="7006  %\Sigma
\mathchardef\zU="7007  %\Upsilon
\mathchardef\zF="7008  %\Phi
\mathchardef\zC="7009  %\Psi
\mathchardef\zW="700A  %\Omega

\newtheorem{theorem}{Theorem}
\newtheorem{corollary}{Corollary}
\newtheorem{proposition}{Proposition}

\newenvironment{pf}{{\noindent{\it Proof. }}}{\ \rule{2mm}{2.5mm}\medskip}

\mathchardef\za="710B  %\alpha
\mathchardef\zb="710C  %\beta
\mathchardef\zg="710D  %\gamma
\mathchardef\zd="710E  %\delta
\mathchardef\zve="710F %\epsilon
\mathchardef\zz="7110  %\zeta
\mathchardef\zh="7111  %\eta
\mathchardef\zvy="7112 %\theta
\mathchardef\zi="7113  %\iota
\mathchardef\zk="7114  %\kappa
\mathchardef\zl="7115  %\lambda
\mathchardef\zm="7116  %\mu
\mathchardef\zn="7117  %\nu
\mathchardef\zx="7118  %\xi
\mathchardef\zp="7119  %\pi
\mathchardef\zr="711A  %\rho
\mathchardef\zs="711B  %\sigma
\mathchardef\zt="711C  %\tau
\mathchardef\zu="711D  %\upsilon
\mathchardef\zvf="711E %\phi
\mathchardef\zq="711F  %\chi
\mathchardef\zc="7120  %\psi
\mathchardef\zw="7121  %\omega
\mathchardef\ze="7122  %\varepsilon
\mathchardef\zy="7123  %\vartheta
\mathchardef\zf="7124  %\varomega
\mathchardef\zvr="7125 %\varrho
\mathchardef\zvs="7126 %\varsigma
\mathchardef\zf="7127  %\varphi
\mathchardef\zG="7000  %\Gamma
\mathchardef\zD="7001  %\Delta
\mathchardef\zY="7002  %\Theta
\mathchardef\zL="7003  %\Lambda
\mathchardef\zX="7004  %\Xi
\mathchardef\zP="7005  %\Pi
\mathchardef\zS="7006  %\Sigma
\mathchardef\zU="7007  %\Upsilon
\mathchardef\zF="7008  %\Phi
\mathchardef\zW="700A  %\Omega

\newcommand{\be}{\begin{equation}}
\newcommand{\ee}{\end{equation}}
\newcommand{\ra}{\rightarrow}

\newcommand{\bea}{\begin{eqnarray}}
\newcommand{\eea}{\end{eqnarray}}
\newcommand{\beas}{\begin{eqnarray*}}
\newcommand{\eeas}{\end{eqnarray*}}
\newcommand{\R}{{\mathbb R}}

\newcommand{\D}{{\rm d}}

\newcommand{\de}{\,{\stackrel{\rm def}{=}}\,}
\newcommand{\we}{\wedge}
\newcommand{\nn}{\nonumber}
\newcommand{\ot}{\otimes}
\newcommand{\s}{{\textstyle *}}

\newcommand{\oX}{\stackrel{o}{X}}
\newcommand{\oD}{\stackrel{o}{D}}
\newcommand{\obD}{\stackrel{o}{\bD}}
\newcommand{\pa}{\partial}
\newcommand{\ti}{\times}
\newcommand{\A}{{\cal A}}

\newcommand{\ka}{\mathbb{K}}

\newcommand{\X}{{\cal X}}

\newcommand{\Ll}{{\pounds}}
\def\lna{\lbrack\! \lbrack}
\def\rna{\rbrack\! \rbrack}

\def\lan{\langle}
\def\ran{\rangle}

\def\ati{{\stackrel{a}{\times}}}
\def\sti{{\stackrel{sv}{\times}}}
\def\aot{{\stackrel{a}{\ot}}}
\def\sati{{\stackrel{sa}{\times}}}
\def\saop{{\stackrel{sa}{\op}}}

\def\svop{{\stackrel{sv}{\oplus}}}
\def\saot{{\stackrel{sa}{\otimes}}}
\def\cti{{\stackrel{cv}{\times}}}
\def\cop{{\stackrel{cv}{\oplus}}}
\def\dra{{\stackrel{\xd}{\ra}}}
\def\bdra{{\stackrel{\bd}{\ra}}}
\def\bAff{\mathbf{Aff}}
\def\Aff{\sss{Aff}}
\def\bHom{\mathbf{Hom}}
\def\Hom{\sss{Hom}}
\def\bt{{\boxtimes}}
\def\sot{{\stackrel{sa}{\ot}}}

\def\op{\oplus}
\def\bwak{{\stackrel{a}{\bigwedge}\!{}^k}}
\def\aop{{\stackrel{a}{\oplus}}}
\def\V{{\cal V}}

\def\cJ{{\cal J}}
\def\bA{\mathbf{A}}
\def\bI{\mathbf{I}}
\def\wh{\widehat}
\def\wt{\widetilde}
\def\ol{\overline}
\def\Sec{\sss{Sec}}
\def\Lin{\sss{Lin}}
\def\ader{\sss{ADer}}
\def\ado{\sss{ADO^1}}
\def\adoo{\sss{ADO^0}}
\def\AS{\sss{AS}}
\def\bAS{\sss{AS}}
\def\bLS{\sss{LS}}
\def\bAP{\sss{AV}}

\def\AP{\sss{AP}}

\def\LS{\sss{LS}}
\def\Z{\mathbf{Z}}
\def\oZ{\overline{\bZ}}

\def\de{{\cal D}^1}
\def\la{\langle}
\def\ran{\rangle}
\def\bd{{\bi d}}

\def\bD{{\bi D}}
\def\bY{{\bi Y}}

\def\bV{{\bi V}}

\def\bS{{\bi S}}

\def\bF{{\bi F}}

\def\bZ{{\bi Z}}

\def\sA{{\sss A}}
\def\sC{{\sss C}}
\def\sD{{\sss D}}

\def\sL{{\sss L}}

\def\sP{{\sss P}}

\def\sT{{\sss T}}
\def\sV{{\sss V}}
\def\sR{{\sss R}}

\def\sF{{\sss F}}

\def\sv{{\sss v}}

\def\xc{\tx{c}}
\def\xd{\tx{d}}
\def\xi{\tx{i}}

\begin{document}
%%%%%%%%%%%%%%%%%%%%%%%%%%%%%%%%%%%%%%%%%%%%%%%%%%%%
\title{AV-differential geometry: Poisson and Jacobi structures\thanks{Research
supported by the Polish Ministry of Scientific Research and
Information Technology under the grant No. 2 P03A 036 25}}
%\footnote{Supported by KBN, grant No 2 P03A 041 18.}
        \author{
        Katarzyna  Grabowska \\
        Division of Mathematical Methods in Physics \\
                University of Warsaw \\
                Ho\.za 69, 00-681 Warszawa, Poland \\
                 {\tt konieczn@fuw.edu.pl} \\
                        \\
        Janusz Grabowski \\
        Institute of Mathematics \\
                Polish Academy of Sciences \\
                ul. \'Sniadeckich 8, P.O.Box 21, 00-956 Warszawa 10, Poland\\
                {\tt jagrab@impan.gov.pl} \\
                        \\
               Pawe\l\ Urba\'nski \\
                Division of Mathematical Methods in Physics \\
                University of Warsaw \\
                Ho\.za 69, 00-681 Warszawa, Poland \\
                {\tt urbanski@fuw.edu.pl}}
\date{}
\maketitle
\begin{abstract}
Based on ideas of W.~M.~Tulczyjew, a geometric framework for a
frame-independent formulation of different problems in analytical
mechanics is developed. In this approach affine bundles replace vector
bundles of the standard description and functions are replaced by sections
of certain affine line bundles called AV-bundles. Categorial constructions
for affine and special affine bundles as well as natural analogs of Lie
algebroid structures on affine bundles (Lie affgebroids) are investigated.
One discovers certain Lie algebroids and Lie affgebroids canonically
associated with an AV-bundle which are closely related to affine analogs
of Poisson and Jacobi structures. Homology and cohomology of the latter
are canonically defined. The developed concepts are applied in solving
some problems of frame-independent geometric description of mechanical
systems.

\bigskip\noindent
\textit{MSC 2000: 17B66 53D10 53D17 70G45}

\medskip\noindent
\textit{Key words: affine spaces; vector bundles; Lie algebroids; Jacobi
structures; Poisson structures; Hamilton formalism.}
\end{abstract}

\section{Introduction} While there is no doubts about the role of
analytical mechanics in explaining many of problems in a variety of
physical topics, it is worth stressing that classical mechanics is by no
means {\em pass\'e}. It is still an open theory with several challenges
and with an influence on both: physics and mathematics. The standard
formulation of analytical mechanics in the language of differential
geometry is based on geometrical objects of vector character. The vector
bundle $\sT M$ of tangent vectors is used as a space of infinitesimal
(dynamical) configurations, the vector bundle $\sT^\s M$ of covectors
plays the role of a phase space, and the Poisson bracket derived from the
symplectic form serve in the Hamiltonian formulation of dynamics in which
one uses the vector space (actually an algebra) of functions. However,
there are situations where one finds difficulties while working with
vector-like objects. Here we list some examples.

\begin{enumerate}
\item  As the first example we describe the problems in the
relativistic mechanics of a charged particle in the external
electromagnetic field. The standard lagrangian $L$ is a function
on the space of infinitesimal configurations $\sT M$: $L(v) = -\la
eA,v\ran + m\sqrt{g(v,v)}$, where $A$ is the one-form representing
the electromagnetic potential, $m$ is the mass and $e$ is the
charge of the particle. The lagrangian depends on the gauge of the
first type. An electromagnetic potential is a connection in the
principal bundle with the structure group $(\R, +)$ over the
space-time. To obtain the one-form representing the potential one
has to choose a section of the bundle (gauge). Changes in the
gauge lead to changes in the lagrangian. The gauge independent
description is possible only when we use affine objects.

\item  The configuration space (the space of events) for the
inhomogenous formulation of time-dependent mechanics is the
space-time $M$ fibrated over the time $\R$. First-jets of this
fibration form the infinitesimal configuration space. Since there
is the distinguished vector field $\pa_t$ on $\R$, the first-jets
of the fibration over time can be identified with those vectors
tangent to $M$ which project on $\pa_t$. Such vectors form an
affine subbundle of the tangent bundle $\sT M$. The bundle $\sV^\s
M$, dual to the bundle of vectors which are vertical with respect
to the fibration over time, is the phase space for the problem.
The phase space carries a canonical Poisson structure, but
hamiltonian fields for this structure are vertical with respect to
the projection on time, so they cannot describe the dynamics. In
the standard formulation the distinguished vector field $\pa_t$ is
added to the hamiltonian vector field to obtain the dynamics. This
can be done correctly when the fibration over time is trivial,
i.e. when $M=Q\ti\R$. When the fibration is not trivial one has to
choose a reference vector field that projects onto $\pa_t$.
Changing the reference vector field means changing the
hamiltonian. To have the description of the dynamic being
independent on the reference field one has to use affine objects.

\item  Let us look on energy and momentum in the most classical
case of Newtonian mechanics. The Newtonian space-time is a
four-dimensional affine space $N$ with an absolute time 1-form
$\zt\in(\sV(N))^\s$, which is a linear function on the model
vector space $\sV(N)$. The dynamics is usually described in a
fixed inertial frame. The inertial frames are represented by
vectors $u\in\sV(N)$ such that $\zt(u)=1$, i.e. $u$ is the
space-time velocity of an inertial observer associated with the
inertial frame. The space of infinitesimal configurations, i.e.
positions and velocities, is $N\ti E_1$, where $E_1\subset\sV(N)$
consists of vectors $v$ satisfying $\zt(v)=1$. Fixing an inertial
frame $u$ allows us to identify $E_1$ with $E_0=Ker(\zt)$ which is
a vector subspace of $\sV(N)$. Therefore we can define momenta as
elements of $E_0^\s=(\sV(N))^\s/\la\zt\ran$. The momentum
transforms according to the formula $p' = p + f(u,u')$ while
changing the inertial frame. The transformation of energy is also
affine, so we cannot describe the dynamics in the
frame-independent way as long as we keep representing the momentum
as a vector object. We need an affine object to replace the usual
covector. We can say that the covector in this case carries too
much structure and we need additional physical information (i.e.
an inertial frame) to use it properly. But even in a fixed
inertial frame the standard description is not satisfactory,
because the identification of $E_1$ with $E_0$ at the very
beginning leads to the use of a wrong Poisson structure to
generate equations of motion from the hamiltonian. This is a
situation similar to the previous example (cf. \cite{GU,Ko}).

\end{enumerate}
Of course, the above list of problems is not complete. Our aim is to
develop the geometric framework for correct approaches. The standard
geometric constructions based on the algebra of functions on a manifold
$M$ are replaced by constructions   based on the affine space of sections
of an affine bundle $\zz\colon \Z\rightarrow M$, modelled on the trivial
bundle $M\times \R$. Such an affine bundle we will call a \textit{bundle
of affine values} (\textit{AV-bundle} in short). The elements of the
bundle $\Z$ replace number-values of functions but we are not informed now
what and where is zero for these values, so our "functions" do not form
any algebra or even a vector space. Such an approach forces deep changes
in the language, notions and canonical objects of differential geometry.
We propose to call this kind of geometry \textit{the differential geometry
of affine values} (\textit{AV-differential geometry} in short).

An additional motivation comes from the observation that even canonical
objects in the traditional "vector geometry" happen to have an affine
character, more or less hidden or forgotten. Let us consider the canonical
symplectic form on the cotangent bundle $\sT^\s M$. This 2-form is
recognized as a linear object while, on the other hand, it is invariant
with respect to translations by closed forms on M that suggests its hidden
affine character. Indeed, it is possible to construct an affine analog of
$\sT^\s M$, which is a symplectic manifold with canonical symplectic
structure and which seems to be more appropriate phase space for many
mechanical problems.

The idea of using affine bundles for the correct frame-independent
geometric formulation of analytical mechanics theories  goes back
to some concepts of W.~M.~Tulczyjew \cite{Tu1} (see also
\cite{Be,TU,Ur1}). We will also use in the paper some of
unpublished ideas of W.~M.~Tulczyjew. A similar approach to
time-dependent non-relativistic mechanics (in the Lagrange
formulation) has been recently developed by E.~Massa with
collaborators \cite{MPL,MVB,Vi}, E.~Mart\'\i nez, W.~Sarlet, and
T.~Mestdag \cite{MMS1,MMS2,MMS3}. Our paper is organized as
follows.

In section 2 we present basic notions of the theory of affine spaces and
relations to the theories of special (resp., cospecial) vector spaces,
i.e. the vector spaces with distinguished a non-zero vector (resp.,
covector).

Basic categorial construction for affine spaces and special/cospecial
vector spaces, like direct sums, products, and tensor products, are
presented in section 3. To our surprise, we could not find them in the
literature.

In section 4, the main affine objects of our approach, namely
\textit{special affine spaces}, i.e. affine spaces modelled on special
vector spaces are introduced together with the corresponding notion of
special duality.

One-dimensional special affine spaces, called \textit{spaces of affine
scalars} are of particular interest. Some properties of such spaces are
investigated in section 5.

All above is extended to the case of bundles in sections 6 and 7. A
\textit{special affine bundle} is a pair $\bA=(A,v)$, where $A$ is an
affine bundle over $M$ and $v\in\Sec(\sV(A))$ is a nowhere-vanishing
section of its model bundle $\sV(A)$. The \textit{dual special affine
bundle} $\bA^\#$ is the affine bundle $\bAff(\bA;\bI)$ of special affine
morphisms of $\bA_m$ into the canonical special affine bundle $M\ti \bI$,
where $\bI=(\R,1)$, i.e. those affine morphisms $\zf:A_m\ra\R$ whose
linear part maps $v(m)$ into $1$, $m\in M$, with the distinguished section
of the model vector bundle being $1_A$ - the constant 1 function on $A$.

One-dimensional special affine bundles are called
\textit{AV-bundles}. An important observation is that there is a
one-to-one correspondence between the space $\Sec(\bA)$ of
sections of a special affine bundle $\bA$ and the space
$\Aff\Sec(\bA^\#)$ of affine sections of the bundle
$\bA^\#\ra\bA^\#/\la 1_A\ran$ which is canonically an AV-bundle.
The affine sections are, of course, those sections $\zs:\bA^\#/\la
1_A\ran\ra\bA^\#$ which are affine maps, i.a. morphisms of affine
bundles. This is a special affine analog of the well-known
correspondence between sections of a vector bundle $E$ and linear
functions on the dual bundle $E^\s$.

In section 8 the phase $\sP\Z$ and the contact bundle $\sC\Z$
associated with an AV-bundle $\Z$ are constructed. They are
AV-analogs of $\sT^\s M$ and $\sT^\s M\op\R$ and carry canonical
symplectic and contact structures, respectively. The AV-Liouville
1-form which is the potential of the canonical symplectic form on
$\sP\Z$ is naturally understood as a section of an affine
fibration over $\sP\Z$ (cf. \cite{Ur2}).

Various Lie algebroids and Lie affgebroids (i.e. Lie algebroid-like
objects on affine bundles \cite{GGU}) associated with a given AV-bundle
$\Z$ are defined and studied in sections 9, 10, 11. Let us mention the Lie
algebroid  $\wt{\sT}\Z$ (an AV-analog of the Lie algebroid extension $\sL
M =\sT M\op\R$ of the canonical Lie algebroid $\sT M$ of vector fields),
the Lie algebroid $\wt{\sL}\Z$ (an AV-analog of the Lie algebroid
extension $\sL M\op\R$ of the Lie algebroid $\sL M=\sT M\op\R$ of linear
first-order differential operators on $M$) and their affine counterparts
$\ol{\sT}\Z$ and $\ol{\sL}\Z$. One proves that the Lie algebroid
$\wt{\sL}\Z$ admits a canonical closed 1-form $\zvf^0$, i.e. $\wt{\sL}\Z$
carries a canonical structure of a \textit{Jacobi algebroid} (see
\cite{IM,GM1,GM2}). It is also shown that sections of $\wt{\sT}\Z$, or
$\wt{\sL}\Z$, (resp., sections of $\ol{\sT}\Z$, or $\ol{\sL}\Z$) can be
interpreted as affine derivations, or affine first-order differential
operators, on sections of $\Z$ with values in functions on $M$ (resp.,
such derivations, or first-order differential operators, but with values
in sections on $\Z$).

In section 12 we recall the definitions and basic facts on
aff-Poisson and aff-Jacobi brackets (cf. \cite{GGU}), i.e. analogs
of Poisson and Jacobi brackets, defined on sections of an
AV-bundle $\Z$ over $M$ and taking values in the ring of smooth
functions $C^\infty(M)$. The main result is the correspondence
between Lie affgebroid structures on a special affine bundle
$\bA=(A,v)$ and aff-Jacobi brackets on the AV-bundle
$\bA^\#\ra\bA^\#/\la 1_A\ran$ which are affine in the sense that
the bracket of two affine sections is an affine function on
$\bA^\#/\la 1_A\ran$. This can be viewed as an AV-analog of the
fact that Lie algebroid brackets on a vector bundle $E$ correspond
to linear Poisson brackets on the dual bundle $E^\s$. In this
picture, the Lagrange formulation of a mechanical problem takes
place on a special affine bundle $\bA=(A,v)$ equipped with a Lie
affgebroid structure, and the lagrangians are sections of the
AV-bundle $\bA\ra\bA/\la v\ran$. The Hamilton formalism, in turn,
takes place on the dual special affine bundle $\bA^\#$ and the
hamiltonians are sections of the AV-bundle $\bA^\#\ra\bA^\#/\la
1_\A\ran$ which carries a canonical aff-Jacobi structure. In most
important examples this structure happens to be aff-Poisson.

In section 13 we observe that aff-Poisson and aff-Jacobi structures on an
AV-bundle $\Z$ correspond to \textit{canonical structures} $\zL$ and $\cJ$
for the Lie algebroid $\wt{\sT}\Z$ and Jacobi algebroid $\wt{\sL}\Z$,
respectively, i.e. $\zL\in\bigwedge^2\wt{\sT}\Z$,
$\lna\zL,\zL\rna_{\wt{\sT}\Z}=0$ (resp., $\cJ\in\bigwedge^2\wt{\sL}\Z$,
$\lna\cJ,\cJ\rna_{\wt{\sL}\Z}^{\zvf^0}=0$), where
$\lna\cdot,\cdot\rna_{\wt{\sT}\Z}$ is the Lie algebroid Schouten bracket
on $\bigwedge^\bullet\wt{\sT}\Z$ (resp.,
$\lna\cdot,\cdot\rna_{\wt{\sL}\Z}^{\zvf^0}$ is the Schouten-Jacobi bracket
of the Jacobi algebroid $(\wt{\sL}\Z,\zvf^0)$). This is an AV-analog of
the well-known identification of Poisson brackets on $C^\infty(M)$ with
Poisson tensors on $M$, i.e. bivector fields with the Schouten-Nijenhuis
square being $0$. The known results on characterization of canonical
structures for Lie and Jacobi algebroids \cite{GM2,GU3} alow one to derive
an analogous characterization for aff-Poisson and aff-Jacobi brackets. In
particular, one can define the corresponding homology and cohomology in a
natural way.

In Section 14 we present solutions of the mentioned problems of
the frame-independent geometric formulation in analytical
mechanics with the use of developed concepts. These solutions form
an alternative to the Kaluza-Klein approach where the vector-like
formulations is kept for the price of extending the dimension (see
also \cite{MT,Tu1,TU}).

Much of this material is to our knowledge new. Our aim was to present a
possibly complete picture which can be viewed as a well-described
mathematical program based on the ideas and needs from analytical
mechanics.

\section{Category of affine spaces}
An {\it affine space} is a triple $(A,V,\za)$, where $A$ is a set,
$V$ is a vector space over a field $\ka$ and $\za$ is a mapping
$\za \colon A \times A\rightarrow V$ such that
      \begin{itemize}
   \item $\za(a_3,a_2) + \za(a_2,a_1) + \za(a_1,a_3) = 0$;
   \item the mapping $\za(\cdot,a) \colon A \rightarrow V$ is bijective for
each $a \in A$.
      \end{itemize}
   We shall also write simply $A$ to denote the affine space $(A,V,\za)$ and
$\sV(A)$ to denote $V$. One can also say that an affine space is a
set with a free  and  transitive action of a vector space (which
is  viewed  as  a  commutative  group  with respect to addition).
By {\it dimension} of $A$ we understand the dimension of $\sV(A)$.
    If $(A,V,\za)$ then also $(A,V,-\za)$ is an affine space. We
will write for brevity $\overline{A}^a$ to denote  the  {\it
adjoint affine space} $(A,V,-\za)$.
   We will write also $a_2 - a_1$ instead of $\za(a_2,a_1)$ and $a + v$ to
denote the unique point $a' \in A$ such that $a' - a = v$,
$v\in\sV(A)$. Of course, every vector  space  is canonically  an
affine space modelled on itself with the affine structure
$\za(v_1,v_2)=v_1-v_2$. The adjoint affine space $\ol{A}^a$ can be
viewed as the same set $A$ with the opposite action of $\sV(A)$:
$a\mapsto a-v$.

It is easy to see that for any linear subspace $V_0$ of $V$ the
set $A/V_0$ of cosets of $A$ with respect to the relation $a\sim
a'\Leftrightarrow  a-a'\in V_0$ is canonically an affine space
modelled on $V/V_0$.

A subset $A'$ of $A$ is an \textit{affine subspace} in $A$ if
there is a linear subspace $\sV(A')$ of $\sV(A)$ such that
$A'=a'+\sV(A')$ for certain $a'\in A'$. Affine subspaces are
canonically affine spaces with the affine structure inherited from
$A$.

If $A'$ is an affine subspace of $A$ then the quotient space $A/A'$ is
understood as $A/\sV(A')$ with distinguished point being the class of $A'$.
Hence $A/A'$ can be identified with the linear space $\sV(A)/\sV(A')$.

Morphisms in the category of affine spaces are affine maps. Let
$A$ and $A'$ be affine spaces. We say that a mapping $\zf\colon
A\rightarrow A'$ is {\it affine} if there is a linear mapping
${\zf}_\sv\colon V\rightarrow V'$ such that
  $$\zf(a + v) = \zf(a) + {\zf}_\sv(v).$$
  We say that ${\zf}_\sv$ is the {\it linear part} of $\zf$.

More  generally, on  every  affine  space instead of the subtraction
$a_1-a_2$, one  can consider {\it vector combination} of elements of $A$,
i.e. the combination $\sum_i\zl_ia_i$, where $a_i\in A$, $\zl_i\in\ka$,
and $\sum_i\zl_i=0$. Every vector combination of elements of $A$ defines a
unique  element  of $\sV(A)$ in obvious way. Similarly, one can consider
{\it affine combinations} (called also {\it barycentric combinations}) of
elements of $A$ which have formally the same form but with
$\sum_i\zl_i=1$. An affine combination determines uniquely an element  of
$A$. Affine maps may be equivalently defined as those maps which respect
affine combinations. Note however that affine combinations do not
determine the affine structure completely: $A$ and $\overline{A}^a$ have
the same affine combinations. The set $\sss{Aff}(A;A')$ of all affine maps
from $A$ into $A'$ is again an affine space modelled on the vector space
$\Aff(A;\sV(A'))$ of affine maps from $A$ into the model vector space
$\sV(A')$ of $A'$: for $\zf_1,\zf_2\in\Aff(A_1;A_2)$ we put
$(\zf_1-\zf_2)(a)=\zf_1(a)-\zf_2(a)$. Inductively, the set
$\sss{Aff}^k(A_1,\dots,A_k;A)$ of $k$-affine maps from $A_1\ti\dots\ti
A_k$ into $A$ is defined as the set
$\sss{Aff}(A_1;\sss{Aff}^{k-1}(A_2,\dots,A_k;A))$. Like in the linear
case, one proves that  $\sss{Aff}^k(A_1,\dots,A_k;A)$ can be identified
with the space of of maps $F:A_1\ti\dots\ti A_k\ra A$ which are affine
with respect to every variable separately. By
$$F^i_\sv:A_1\ti\dots\ti\sV(A_i)\ti\dots\ti A_k\ra\sV(A)
$$
we denote the linear part of $F$ with respect to the $i$th
variable. It is linear on $\sV(A_i)$ and affine with respect to
each of the remaining variables separately. The higher-order
linear parts $F^{i_1,\dots,i_l}_\sv$ are defined in obvious way.
The multilinear map
$$F^{1,\dots,k}_\sv:\sV(A_1)\ti\dots\ti\sV(A_k)\ra\sV(A)
$$
we denote simply by $F_v$.

A  {\it  free  affine  space} $\A=\A(\{ a_j\}_{j\in J})$ generated
by the set $\{ a_j:j\in J\}$ is an affine subspace in the free
vector space generated by $\{ a_j:j\in J\}$, i.e. in the the
vector space $\V(\{  a_j\}_{j\in J})$ of formal linear
combinations $v=\sum_j\zl_ja_j$,    described    by the equation
$1_\A(v)=1$, where $1_\A$ is the linear functional on $\V(\{
a_j\}_{j\in J})$ defined by  $1_\A(v)=\sum_j\zl_j$. The notation
is justified by the fact that this functional is constantly 1 on
$\A$. The model vector space for this free affine space is a
linear subspace of $\V(\{ a_j: j\in J\})$ being the kernel of the
functional $1_\A$. Every affine space $A$ is actually isomorphic
to the free affine space generated by a subset of $A$ which we
call a {\it basis} of $A$. A subset $\{ a_j:j\in J\}$ is a basis
of $A$ if every element of $A$ can be expressed uniquely as an
affine combination of elements of the basis. Existence of a basis
can be proved analogously to the linear case, since  $\{ a_j:j\in
J\}$ is a basis of $A$ if and only if  $\{ a_j-a_{j_0}:j\in J'\}$
is a basis  of $\sV(A)$, where $j_0\in J$ and $J'=J\setminus \{
j_0\}$. The dimension of $A$ is the cardinality of a basis minus
1.

Every affine space $A$ is canonically embedded as an affine hyperspace into a
vector space $\widehat{A}$ which we call the {\it vector hull} of $A$. The
vector  hull  $\widehat{A}$  is defined  as  the quotient space
$\V(A)/\V_0(A)$ of the free vector space $\V(A)$ generated by  $A$ by  its
subspace spanned by    linear    combinations     of the form
$1\cdot(a+\zl(a'-a''))-1\cdot a-\zl a'+\zl a''$. Here the expression
$(a+\zl(a'-a''))$ is viewed as an element of $A$. Since $A$ is canonically
embedded into $\V(A)$ as a set, we have a canonical  map from $A$ into
$\widehat{A}$ which can be proved to be an embedding of the affine space onto
an {\it affine hyperspace}, i.e. a 1-codimensional affine subspace which is
proper (does not contain 0), of $\widehat{A}$. This hyperspace can be
equivalently defined as   the level-1 set of the functional
$1_A:\widehat{A}\ra\ka$ represented  by the sum of coefficients on $\V(A)$.
We will not denote this embedding in a special way just regarding $A$ as a
subset of $\wh{A}$. The model vector space $\sV(A)$ is also canonically
embedded in $\wh{A}$ as the kernel of $1_A$.

Choosing a basis $\{ a_j:j\in J\}$ of $A$ we get an isomorphism of
$\widehat{A}$ with $\V(\{ a_j:j\in J\})$. Note that for a vector space $V$
viewed as an affine space its vector hull $\widehat{V}$ is canonically
isomorphic to $V\oplus\ka$. This decomposition follows from the existence
of a distinguished element $0\in V$ which is a non-zero vector in $\wh{V}$
complementary to $\sV(V)\simeq V$. It is obvious by construction that the
vector hull is unique up to isomorphism, so that we have the following.
\begin{theorem} Every affine space $A$ is canonically embedded as
an affine hyperspace of the vector space $\wh{A}$ -- its vector hull.
Conversely, if $A$ is embedded as an affine hyperspace of a vector space
$W$, then there is a canonical isomorphism $\zF:\wh{A}\ra W$ which reduces
to the identity map on the embedded $A$.
\end{theorem}
For vector spaces $V_1,V_2$ we denote by $\Hom(V_1;V_2)$ the space
of morphisms (linear maps) from $V_1$ into $V_2$ and by
$\Hom^{A_2}_{A_1}(V_1;V_2)$ the subset of those morphisms
$\zF\in\Hom(V_1;V_2)$ which map the subset $A_1$ of $V_1$ into the
subset $A_2$ of $V_2$.
\begin{theorem}\label{al} For an affine space $A$ and a vector space $V$
there are canonical identifications
\begin{description}
\item{(a)}
$$\Aff(A,V)\ni\zf\mapsto\widehat{\zf}\in\Hom(\widehat{A},V).$$
In particular, the vector space  $A^\dag=\Aff(A,\R)$ is
canonically isomorphic to $\wh{A}^*$, and
\item{(b)}
$$\Aff(A_1,A_2)\ni\zf\mapsto\widehat{\zf}\in\Hom^{A_2}_{A_1}
(\widehat{A_1},\wh{A_2})$$ for affine spaces $A_1,A_2$.
\end{description}
\end{theorem}

\begin{pf}
We put simply $\widehat{\zf}(\sum_i\zl_ia_i)=\sum_i\zl_i\zf(a_i)$ for
$\zl_i\in\ka$, $a_i\in A$, so $\zf=\widehat{\zf}_{\mid A}$. There are
obvious embeddings
$$\Aff(A_1,A_2)\subset\Aff(A_1,\wh{A_2})\subset\Hom(\wh{A_1},\wh{A_2})$$
and it is easy to see that $\Aff(A_1,A_2)$ is characterized inside
$\Hom(\wh{A_1},\wh{A_2})$ as the set of those morphisms which map $A_1$
into $A_2$.
\end{pf}

\medskip
The vector space $\widehat{A}$ has a distinguished affine
hyperspace $A$. Such an affine subspace is uniquely determined by
the nonzero functional  $1_A\in\widehat{A}^*$ as its level-1 set:
$A=\{ v\in\widehat{A}:1_A(v)=1\}$. Thus
$\widehat{\bA}=(\widehat{A},1_A)$ is an example of a {\it
cospecial vector space}, i.e. a vector space  with a distinguished
affine hyperspace, or, equivalently, as a vector space with a
distinguished non-zero linear functional. On the other hand, its
vector dual $\widehat{A}^*$, which is canonically identified with
$A^\dag=\Aff(A,\ka)$, is a {\it special vector space}, i.e. a
vector space with a distinguished non-zero element. We will denote
this special vector space by $\bA^\dag=(A^\dag,1_A)$ and call it
the {\it vector dual} of $A$. In finite dimension we have a true
duality between affine spaces and special vector spaces. Indeed,
every special vector space $\bV=(V,v^0)$ defines an affine
hyperspace $\bV^\ddag=\{ u\in V^*:u(v^0)=1\}$ in the dual $V^*$ of
$V$. Since in finite dimension $(V^*)^*=V$, we have the following.
\begin{theorem}
For finite-dimensional special vector space $\bV$ and
finite-dimensional affine space $A$ there are canonical
isomorphisms
$$((\bV^\ddag)^\dag,1_{\bV^\ddag})\simeq\bV
$$
and
$$(\bA^\dag)^\ddag\simeq A.
$$
\end{theorem}

The vector hull $\widehat{\sss{Aff}}(A_1,A_2)$ of $\sss{Aff}(A_1,A_2)$ can
be interpreted as the vector space
$\widehat{\sss{Hom}}(\widehat{\bA_1},\widehat{\bA_2})$ of those linear
maps $F:\widehat{A_1}\ra\widehat{A_2}$ for which $F^*(1_{A_2})=\zl
1_{A_1}$ for certain $\zl\in\ka$.

Special (resp., cospecial) vector spaces form a category with the
set of morphisms $\bHom(\bV_1,\bV_2)$ between $\bV_i=(V_i,v_i^0)$
(resp. $\bV_i=(V_i,\zf_i)$), $i=1,2$, consisting of those linear
maps $F:V_1\ra V_2$ for which $F(v_1^0)=v_2^0$ (resp.,
$F^*(\zf_2)=\zf_1$). The condition $F^*(\zf_2)=\zf_1$ means that
$F$ maps the points of the affine hyperspace
$A_1=\{\zf_1(u_1)=1\}$ of $V_1$ into the affine hyperspace
$A_2=\{\zf_2(u_2)=1\}$ of $V_2$. There is a canonical covariant
equivalence functor from the category of cospecial vector spaces
into the category of affine spaces. It associates with any
cospecial vector space $(V,A)$ its affine hyperspace $A$, and with
every morphism $F:(V_1,A_1)\ra(V_2,A_2)$ its restriction to $A_1$
(which is an affine map into $A_2$). Conversely, with every affine
space $A$ we associate its vector hull $\widehat{A}$ with $A$ as
the distinguished affine hyperspace and with every affine map
$F:A_1\ra A_2$ its (unique) extension to a linear map from
$\widehat{A_1}$ into $\widehat{A_2}$. In finite dimensions we can
use the duality and obtain a contravariant equivalence functor
from the category of special vector spaces to the category of
affine spaces. This functor associates with a special vector space
$\bV=(V,v^0)$ the affine hyperspace $\bV^\ddag$ in $V^*$. We will
use these equivalences to construct categorial object for the
affine category exploring (generally better) knowledge of the
linear category.

\section{Categorial constructions for affine spaces}

In the category of special vector spaces and, consequently, in the
category of cospecial vector spaces and the category of affine
spaces there are direct sums and products. We will just describe
the models leaving the obvious proofs to the reader. The
constructions are very natural but, to our surprise, we could not
find explicit references in the literature.

For special vector spaces $\bV_i=(V_i,v^0_i)$, $i=1,2$, their product
$\bV_1\sti\bV_2$ is represented by the standard product $V_1\ti V_2$ with
the distinguished vector $v^0=(v^0_1,v^0_2)$. The projections
$\zp_i:V_1\ti V_2\ra V_i$ map $v^0$ onto $v^0_i$, i.e. represent morphisms
of special vector spaces.

The {\it special direct} sum $\bV_1\svop\bV_2$ is represented by
the quotient vector space $V_1\oplus V_2/\langle
v^0_1-v^0_2\rangle$ with the distinguished vector being the class
$[v^0_1]$ of $v^0_1$ (or, equivalently, the class $[v^0_2]$ of
$v^0_2$). The embedding of $\bV_i$ is represented by the embedding
of $V_i$ in $V_1\op V_2$ composed with the projection.

By duality, for cospecial vector spaces $\bV_i=(V_i, \zf^0_i)$,
$i=1,2$, its {\it cospecial direct sum} $\bV_1\cop\bV_2$ is
represented by the vector space $V_1\oplus V_2$ with the
distinguished functional $\zf^0=(\zf^0_1,\zf^0_2)\in V_1^*\ti
V_2^*=(V_1\op V_2)^*$ and obvious embeddings of $\bV_i$. The
product $\bV_1\cti\bV_2$, in turn, is represented by the linear
hyperspace in $V_1\ti V_2$ being the kernel of $\zf^0_1-\zf^0_2\in
V_1^*\op V_2^*=(V_1\ti V_2)^*$ and equipped with the distinguished
functional $\zf^0=(\zf^0_1)_{\mid
Ker(\zf^0_1-\zf^0_2)}=(\zf^0_2)_{\mid Ker(\zf^0_1-\zf^0_2)}$. The
projections from $Ker(\zf^0_1-\zf^0_2)$ onto $V_i$ are just
restrictions of projections from $V_1\ti V_2$. They give rise to
cospecial morphisms from $\bV_1\cti\bV_2$ onto $\bV_i$.

The above constructions allow us to recognize the products and
sums in the category of affine spaces. The {\it affine product}
$A_1\ati A_2$ in this category is the standard cartesian product
$A_1\ti A_2$ which is an affine space modelled on
$\sV(A_1)\ti\sV(A_2)$,
$(a_1,a_2)-(a_1',a_2')=(a_1-a_1',a_2-a_2')$. The direct sum
$A_1\aop A_2$ is the affine hyperspace in $\wh{A_1}\op\wh{A_2}$
generated by the affine subspaces $A_1,A_2$ which are canonically
embedded, i.e.
$$A_1\aop A_2=\{\zl_1a_1+\zl_2a_2\in\wh{A_1}\oplus\wh{A_2}: a_1\in A_1,
a_2\in A_2, \zl_1+\zl_2=1\},
$$
with obvious embeddings of $A_1$ and $A_2$.

\begin{theorem}\label{cc} We have canonical isomorphisms
\bea \widehat{(A_1\ati A_2)}&\simeq &\widehat{A_1}\cti\widehat{A_2},\\
\widehat{(A_1\aop A_2)}&\simeq &\widehat{A_1}\cop\widehat{A_2},\\
(A_1\ati A_2)^\dag&\simeq &A_1^\dag\svop A_2^\dag,\\
(A_1\aop A_2)^\dag&\simeq &A_1^\dag\sti A_2^\dag,\eea where the
vector hulls and the vector duals are regarded as cospecial and
special vector spaces, respectively.
\end{theorem}

In the category of affine spaces we can define {\it affine tensor
products} $A_1\aot\dots\aot A_k$ which are affine spaces such that
$\sss{Aff}^k(A_1,\dots,A_k;A)=\sss{Aff}(A_1\aot\cdots\aot A_k;A)$. Like in
the linear case, $A_1\aot\cdots\aot A_k$ can be defined as the quotient of
the free affine space $\A(A_1\ti\dots\ti A_k)$ by the linear subspace of
its model vector space generated by elements of the form
$$(a_1,\dots,a_i+\zl(a_i'-a_i''),\dots,a_k)-
(a_1,\dots,a_i,\dots,a_k)-\zl(a_1,\dots,a'_i,\dots,a_k)
+\zl(a_1,\dots,a''_i,\dots,a_k).$$  One can also say that
$A_1\aot\cdots\aot A_k$ is the affine subspace in
$\widehat{A_1}\ot\dots\ot\widehat{A_k}$ spanned by tensors of the
form $a_1\ot\dots\ot a_k$, where $a_i\in A_i$. The tensor product
$A_1\aot\dots\aot A_k$ may be viewed also as the affine hyperspace
in the standard tensor product
$\widehat{A_1}\ot\dots\ot\widehat{A_k}$ determined by the
functional $1_{A_1}\ot\dots\ot 1_{A_k}$ and the associated vector
space $\sV(A_1\aot\cdots\aot A_k)$ is the kernel of
$1_{A_1}\ot\dots\ot 1_{A_k}$. It is easy to see that
 $\sV(A_1\aot\cdots\aot A_k)$ is additively generated by tensors
$a_1\ot\dots\ot v_i\ot\dots\ot a_k$ from $\widehat{A_1}\ot\dots\ot
\widehat{A_k}$, where $a_j\in A_j$, $v_i\in\sV(A_i)$. This is indeed a
vector space, since $\zl(a_1\ot\dots\ot v_i\ot\dots\ot a_k)$ is
represented by $a_1\ot\dots\ot\zl v_i\ot\dots\ot a_k$. If we fix $a^0_i\in
A_i$, then
$$A_1\aot\dots\aot A_k=a^0_1\ot\dots\ot a^0_k+\bigoplus_{i_1<\dots<i_l}
a^0_1\ot\dots\ot\sV(A_{i_1})\ot\dots\ot\sV(A_{i_l})\ot\dots\ot
a^0_k.
$$
Sometimes we will write $a_1\aot\dots\aot a_k$ for the affine
tensor product represented by $a_1\ot\dots\ot
a_k\in\widehat{A_1}\ot\dots\ot\widehat{A_k}$ to stress that we are
dealing with an element of $A_1\aot\dots\aot A_k$. The canonical
map
$$A_1\ti\dots\ti A_k\ni (a_1,\dots,a_k)\mapsto a_1\aot\dots\aot
a_k\in A_1\aot\dots\aot A_k$$ is a multi-affine map. Note that for
vector spaces $V_i$ there is a canonical identification of
$V_1\aot\dots\aot V_k$ with
$((V_1\oplus\ka)\ot\dots\ot(V_k\oplus\ka))\ominus(\ka\ot\dots\ot\ka)$.
For the dimension we have the formula
$$dim(A_1\aot\cdots\aot A_k)=(dim(A_1)+1)\cdot\dots\cdot(dim(A_k)+1)-1.
$$
Like in the linear case, we have natural isomorphisms

\bea A_1\aot A_2&\simeq& A_2\aot A_1,\\
 A_1\aot(A_2\aot A_3)&\simeq&(A_1\aot A_2)\aot A_3,\\
\widehat{A_1\aot A_2}&\simeq& \widehat{A_1}\ot\widehat{A_2},\\
\label{111}(A_1\aot A_2)^\dag&\simeq& (A_1)^\dag\ot(A_2)^\dag.
\eea To define affine skew-symmetric tensor product $\bwak A$ for
$k>1$, let us observe first that the symmetric group $S_k$ acts
naturally on $A^{\aot k}$. By $A^k_0$ we denote its affine
subspace spanned by tensors of the form $a_1\aot\dots\aot a_k$,
where $a_i=a_j$ for certain $i\ne j$, i.e. invariant with respect
to a transposition. We put
$$\bwak A=A^{\aot k}/A^k_0
$$
which is canonically a vector space (see the previous section). It follows
directly from definition that any element of $\sss{Aff}(\bwak A;A')$
represents a multi-affine mapping $F:A\ti\dots\ti A\ra A'$ (an element of
$\Aff(A^{\aot k};A')$) which is constant on $A^k_0$, i.e. constant on the
set of those $(a_1,\dots,a_k)$ for which $a_i=a_j$ for certain $i\ne j$.

It is a standard task to prove that such multi-affine mappings are {\it
skew-symmetric} in the sense that
$$F^{\zs^{-1}(1)}_\sv((a_{\zs^(1)},\dots,v,\dots,a_{\zs(k)})=sgn(\zs)
F^1_\sv(v,a_2,\dots,a_k)
$$
for any permutation $\zs\in S_k$ and $a_2,\dots,a_k\in A$,
$v\in\sV(A)$. It is easily seen that, for $k>1$, the affine wedge
product $\bwak A$ is canonically isomorphic to the standard
exterior power $\bigwedge^k\widehat{A}$. To put it simpler, one
can also say that $\bwak A$ is the affine subspace in
$\bigwedge^k\widehat{A}$ generated by tensors $a_1\we\dots\we
a_k$, $a_i\in A$, which happens to be the whole
$\bigwedge^k\widehat{A}$.

\section{Special affine spaces and special duality}
A {\it special affine space} $\bA=(A,v^0)$ is an affine space $A$
modelled on a special vector space  $\bV(A) =(\sV(A),v^0)$. The {\it
adjoint special affine space} $\overline{\bA}=({A},-v^0)$ is
modelled on the {\it adjoint special vector space}
$\ol{\bV}=(V,-v^0)$.

Let $\bA=(A,v^0)$ and $\bA_i = (A_i,v^0_i)$, $i=1,\dots,k$, be
special affine spaces with the distinguished vectors
$v^0\in\bV(A)$, $v_i^0\in \bV(A_i)$. By $\bAff(\bA_1;\bA)$  we
denote the space of special affine maps $\zf:\bA_1\ra\bA$. It is
canonically a special affine space, since the constant map onto
$\{ v^0\}$ is naturally distinguished in
$\sV(\Aff(A_1;A))=\Aff(A_1,\sV(A))$.
    Inductively, we put
$$\bAff^k(\bA_1,\dots,\bA_k;\bA)=\bAff(\bA_1;\bAff^{k-1}(\bA_2,
\dots,\bA_k;\bA))
$$
for the space of $k$-special affine  maps  from $\bA_1\ti\dots\ti
\bA_k$ into $\bA$, which consists of maps $F:A_1\ti\dots\ti A_k\ra
A$ which are special affine with respect to every variable
separately.

The vector hull of a special affine space is canonically a {\it bispecial
vector space}, i.e. a vector space $V$ with a distinguished non-zero vector
$v^0\in V$ and a distinguished non-zero covector $\zf^0\in V^*$ (or an affine
hyperspace $A$) such that $\zf^0(v^0)=0$ (or $v^0\in\sV(A)$). Morphisms
between bispecial vector spaces $\bV_i=(V_i,v^0_i,\zf^0_i)$, $i=1,2$, are those
linear maps $F:V_1\ra V_2$ which respect the distinguished vectors and
covectors: $F(v^0_1)=v^0_2$, $F^*(\zf^0_2)=\zf^0_1$.

In finite dimensions we have the obvious equivalence between the
category of special affine spaces and the category of bispecial
vector spaces. Since the category of bispecial vector spaces,
which is canonically equivalent to the category of special affine
spaces, is self-dual with the obvious duality
$(V,v^0,\zf^0)^\#=(V^*,\zf^0,v^0)$, we have the natural duality
$\bA\leftrightarrow\bA^\#$ in the category of special affine
spaces. The dual $\bA^\#$ of the special affine space $\bA=(A,v^0)$
is thus the affine hyperspace in $(\widehat{A})^*=A^\dag$ defined
as the level-$1$ set of the functional $v^0$ (we use the embedding
$\widehat{A}\subset\widehat{A}^{**}$) and equipped with the vector
$1_A$, $1_A(A)=1$, of its model vector space. In other words,

$$ \bA^\#=\bAff(\bA,\bI),$$ where $\bI=(\ka,1)$ is the canonical
special vector space, with canonically chosen map $1_A$ in
$\sV(\bAff(\bA,\bI))$. Let us observe that $1_A$ really belongs to
the model vector space for $\bA^\#$, since the latter consists of
those affine maps $\zf:A\ra\ka$ whose linear part vanishes on
$v^0$, i.e.
\beas \sV(\bA^\#)&=&\{\zf\in\Aff(A,\ka): \zf_\sv(v^0)=0\}=
\Aff(A/\langle v^0\rangle;\ka)=\\ &&\{\wh{\zf}\in\Hom(\wh{A};\ka):
\wh{\zf}(v^0)=0\}= \Hom(\wh{A}/\langle v^0\rangle;\ka).
\eeas
The bispecial interpretation of the special affine duality yields
immediately the canonical isomorphism $(\bA^\#)^\#=\bA$. Note also
that one can view $\bI$ as $\{ *\}^\dag$, where $\{ *\}$ is a
single-point affine space, and that the map $\zf\mapsto -\zf$
establishes a canonical isomorphism
$$\ol{\bA^\#}\simeq\ol{\bA}^\#.$$

A {\it special affine pairing} between special affine spaces
$\bA_1$ and $\bA_2$ is a special biaffine map
$$\zF:\bA_1\ti\bA_2\ra\bI$$
for which the corresponding maps
$$\zF^l:\bA_1\ra\bA_2^\#=\bAff(\bA_2,\bI),\quad
\zF^l_{a_1}(a_2)=\zF(a_1,a_2),$$ and
$$\zF^r:\bA_2\ra\bA_1^\#=\bAff(\bA_1,\bI),\quad
\zF^r_{a_2}(a_1)=\zF(a_1,a_2)$$ are isomorphisms (in finite dimension it
is sufficient that they are injective). An example is given by the
canonical special affine pairing of dual special affine spaces

$$ \langle\cdot,\cdot\rangle_{sa}:\bA\ti\bA^\#\ra\bI,\quad \langle
a,\zf\rangle_{sa}=\zf(a)=a(\zf)
$$
This is just the restriction of the pairing between $\wh{\bA}$ and
$\wh{\bA^\#}=\bA^\dag=\wh{\bA}^*$ to the product of affine
hyperspaces $\bA\ti\bA^\#$. Note that every special affine map
$\zc\in\bAff(\bA_1;\bA_2)$ has its dual
$\zc^\#\in\bAff(\bA_2^\#;\bA_1^\#)$ defined by
$$\langle
a_1,\zc^\#(a_2^\#)\rangle_{sa}=\langle\zc(a_1),a_2^\#\rangle_{sa}.
$$
Note also that the concept of special vector spaces and the corresponding
duality has been introduced in \cite{Ur1}.

Since morphisms of bispecial vector spaces are linear maps which are
simultaneously morphisms of special and cospecial structures, we can
combine the constructions of the previous section to get products, direct
sums, and tensor products in the category of special affine spaces.

Recall that the special direct sum $\bV_1\svop\bV_2$ is
represented by the quotient vector space $V_1\oplus V_2/\langle
v^0_1-v^0_2\rangle$ with the distinguished vector being the class
$[v^0_1]$ of $v^0_1$ (or, equivalently, the class $[v^0_2]$ of
$v^0_2$). A similar construction $\bA_1\bt\bA_2=((A_1\ti A_2)/\lan
(v^0_1,-v^0_2)\ran,[(v^0_1,v^0_2)])$ we can perform in the
category of special affine spaces. The model space for
$\bA_1\bt\bA_2$, which will be called the {\it reduced product},
is canonically isomorphic with $\sV(\bA_1)\svop\sV(\bA_2)$.
However, $\bA_1\bt\bA_2$ is not the direct sum in the category of
special affine spaces which will be constructed in a while. The
class of $u+v$ in $\sV(\bA_1)\svop\sV(\bA_2)$ (resp., the class of
$(a_1,a_2)$ in $\bA_1\bt\bA_2$) we will denote by $u\svop v$
(resp., $a_1\bt a_2$). Note that any special affine pairing
$\zF:\bA_1\ti\bA_2\ra\bI$ is constant on fibers of the canonical
projection $\bA_1\ati\bA_2\ra\bA_1\bt\bA_2$. The notion of the
reduced product is useful because of the following fact which can
be easily derived from Theorem \ref{cc} (3).

\begin{theorem}\label{bt} For special affine spaces $\bA_i$, $i=1,2$, we
have
$$(\bA_1\bt\bA_2)^\#\simeq \bA_1^\#\bt\bA_2^\#.$$
In particular, for any affine space $A$ and $\bA_1=A\ati\bI$, one
has $A\ati\bA_2=\bA_1\bt\bA_2$ and consequently
$$(A\ati\bA_2)^\#\simeq A^\dag\bt\bA_2^\#.$$
\end{theorem}

For special affine spaces $\bA_i=(A_i,v^0_i)$, $i=1,2$, their {\it
special affine direct sum} $\bA_1\saop\bA_2$ is represented by the
affine space $(A_1\aop A_2)/\langle v^0_1-v^0_2\rangle$ modelled
on $Ker(1_{A_1}+1_{A_2})/\langle v^0_1-v^0_2\rangle$ in
$(\wh{A_1}\op\wh{A_2})/\langle v^0_1-v^0_2\rangle$ with the
distinguished vector $v^0\in(\wh{A_1}\op\wh{A_2})/\langle
v^0_1-v^0_2\rangle$ being the class $[v^0_1]$ of $v^0_1$ (or,
equivalently, the class  $[v^0_2]$ of $v^0_2$). There are obvious
special affine embeddings of $\bA_i$ into $\bA_1\saop\bA_2$,
$i=1,2$.

The {\it special affine product} $\bA_1\sati\bA_2$ is represented
by $A_1\ati A_2$ modelled on $\sV(\bA_1)\sti\sV(\bA_2)$ with
distinguished vector $v^0=(v_1^0,v_2^0)$ in
$\sV(\bA_1)\ti\sV(\bA_2)$. The special affine projections from
$\bA_1\sati\bA_2$ onto $\bA_i$, $i=1,2$, are obvious. Note that
the dimensions of $\bA_1\saop\bA_2$ and $\bA_1\sati\bA_2$ are
equal, but the model vector spaces are different (we have
inclusion $\sV(\bA_1)\ti\sV(\bA_2)\subset Ker(1_{A_1}+1_{A_2})$,
but $v^0_1-v^0_2\in \sV(\bA_1)\ti\sV(\bA_2)$). However, like for
vector spaces, they are related by duality.
\begin{theorem} There are canonical isomorphisms
\bea(\bA_1\sati\bA_2)^\#&\simeq &\bA_1^\#\saop\bA_2^\#,\\
(\bA_1\saop\bA_2)^\#&\simeq &\bA_1^\#\sati\bA_2^\#. \eea
\end{theorem}

For special multi-affine morphisms from $\bA_1\ti\dots\ti\bA_k$ we have a
representing object, the {\it special affine tensor product}
$\bA_1\saot\dots\saot\bA_k$, such that
$$\bAff^k(\bA_1,\dots,\bA_k;\bA)=\bAff(\bA_1\saot\dots\saot\bA_k;\bA).$$
This is the quotient of the affine tensor product
$A_1\aot\dots\aot A_k$ by the linear subspace of
$\sV(A_1\aot\cdots\aot A_k)$ spanned by tensors
$$a_1\aot\dots\aot v^0_i\aot\dots\aot a_k-b_1\aot\dots\aot v^0_j\aot\dots\aot b_k,
$$
where $a_l,b_l\in A_l$.  In the special affine case
$dim(\bA_1\sot\bA_2)=dim(\bA_1)\cdot dim(\bA_2)$. The canonical
map
$$\bA_1\ti\dots\ti\bA_k\ni (a_1,\dots,a_k)\mapsto a_1\saot\dots\saot
a_k\in \bA_1\saot\dots\saot\bA_k,
$$
where $a_1\saot\dots\saot a_k$ is the coset of $a_1\aot\dots\aot a_k$, is
a special multi-affine map. We have obvious canonical isomorphisms

$$\bA_1\saot\bA_2\simeq\bA_2\saot\bA_1, \quad
\bA_1\saot(\bA_2\saot\bA_3)\simeq(\bA_1\saot\bA_2)\saot\bA_3,
\quad (\bA_1\saot\bA_2)^\#=\bA_1^\#\saot\bA_2^\#.
$$

Note that there are no special $k$-affine and skew-symmetric maps
$F:\bA\ti\dots\ti\bA\ra\bA'$, $\bA=(A,v^0)$, $\bA'=(A',v')$ for $k>1$, since
special $k$-affine implies that
$F^1_\sv(v^0,a,\dots)=F^2_\sv(a,v^0,\dots)=v'$ and skew-symmetry that
$F^1_\sv(v^0,a,\dots)=-F^2_\sv(a,v^0,\dots)$, so $v'=0$ ; a contradiction.

Starting with an object in vector or affine category we can always
construct canonically an object in the special category just by taking the
product with $\bI=(\R,1)$. For example, given an affine space $A$ we can
define its {\it specialization} $\bS_A$ as the special affine space
$\bS_A=(A\ati\bI, (0,1))$ modelled on the specialization
$\bS_{\sV(A)}=(\sV(A)\ti\bI, (0,1))$ of the model space for $A$. Using the
specialization we can describe certain canonical constructions in affine
category in the language of the special affine category. Note that in this
language $\bI=\bS_{\{ *\}}$, where $\{ *\}$ is a one-point affine space.

\begin{theorem} For affine spaces $A,A_1,A_2$ there are canonical
isomorphisms
\bea(a)\qquad\qquad \bS_A^\#&\simeq &A^\dag,\\
(b)\qquad\qquad \wh{\bS_A}&\simeq &\bS_{\wh{A}},\\
(c)\qquad \bS_{A_1\aop A_2}&\simeq &\bS_{A_1}\saop\bS_{A_2},\\
(d)\qquad \bS_{A_1\ati A_2}&\simeq &\bS_{A_1}\bt\bS_{A_2}.\eea
\end{theorem}
\begin{pf} The proof is straightforward and we will prove only (a)
leaving the rest to the reader. Since any
$\zf\in\bS_A^\#=\bAff(\bS_A;\bI)$ is an affine map characterized
by $\zf(a,r)=\zf(a,0)+r$, there is a one-one correspondence
$$\bS_A^\#\ni\zf\mapsto\zf^\dag\in A^\dag$$ given by
$\zf^\dag(a)=\zf(a,0)$. It is obvious that
$(\zf+1_{\bS_A})^\dag=\zf^\dag+1_A$, so this is an isomorphism of
special affine spaces.
\end{pf}

\section{Spaces of affine scalars}
In this section we will consider 1-dimensional special affine
spaces $\Z=(Z,v^0)$. Since the model vector space $\sV(\Z)$ is
1-dimensional and special, we can identify it with $\bI=(\ka,1)$.
In what follows we assume that $\Z$ is modelled on $\bI$. In this
picture the adjoint special affine space $\ol{\Z}$ is isomorphic
to $(\ol{Z}^a,1)$, i.e. $\ol{\Z}$ is $Z$ with the same
distinguished vector but with the adjoint affine structure:
$\zs-_0\zs'=\zs'-\zs$ (or $\zs+_0r=\zs+(-r)$). The points of $\Z$
are like numbers, i.e. elements of $\ka$, but the origin 0 is not
fixed, so only the difference of points makes sense as a number
or, equivalently, we can add numbers to points of $\Z$. We will
call $\Z$ a {\it space of affine scalars}. Of course, any point
$\zs_0\in\Z$ defines the isomorphism $I_{\zs_0}:\Z\ra\bI$,
$\zs\mapsto \zs-\zs_0$, of special affine bundles. We can consider
also the map $\bF:Z\ra Z^\dag=\Aff(Z,\ka)$ given by
$\bF_\zs(\zs')=\zs-\zs'$. The following is straightforward.

\begin{theorem}\label{tf} The map $\zs\mapsto\bF_\zs$
induces a canonical isomorphism
$\bF:\Z\ra\oZ^\#$ represented by the special affine pairing
$$\Z\ti\oZ\ni(\zs,\zs')=\bF_\zs(\zs')=\zs-\zs'\in\bI.
$$
This isomorphism extends by linearity to an isomorphism
$\bF:\wh{\Z}\ra\Z^\dag$ of special vector spaces.
\end{theorem}

There are canonical geometric structures on the space of affine
numbers $\Z$. Since a translation of a polynomial function on
$\ka$ is a polynomial function, the algebra $Pol(\Z)$ of
polynomial functions on $\Z$ is well-defined. It is generated by
affine functions on $\Z$. There is a canonical `vector field'
(derivation of $Pol(\Z)$) on $\Z$ being the `fundamental vector
field' $X_\Z$ of the $\ka$-action on $\Z$, $\zs\mapsto\zs+s$. With
respect to any `global coordinate system' $I_\zs:\Z\ra\ka$ this
vector field has the form $X_\Z=-\partial_s$, where
$\partial_s(s^n)=ns^{n-1}$ for $s$ being the standard coordinate
in $\ka$. We can also consider a Jacobi structure on $\Z$ with the
corresponding Jacobi bracket
\be\label{jb} \{ f,g\}_\Z=fX_\Z(g)-gX_\Z(f)
\ee on $Pol(\Z)$. Of course, in the case $\ka=\R$ one can
understand $X_\Z$ as a true vector field on $\Z$ and the bracket
$\{\cdot,\cdot\}_\Z$ can be understood as a bracket defined on the
algebra $C^\infty(\Z)$ of all smooth functions on $\Z$.

\begin{proposition} \begin{description} \item{(a)} For all
$\zs,\zs'\in\Z$,
\be\label{jbdu} \{ \bF_\zs,\bF_{\zs'}\}_\Z=\bF_\zs(\zs')=\zs-\zs';\ee
\item{(b)} For all $\zvf\in\Z^\dag$ and all $\zs\in\Z$,
\be\label{jbdu1} \{ \zvf,\bF_\zs\}_\Z=\zvf(\zs).\ee
\end{description}
\end{proposition}
\begin{pf} Let us identify $\Z$ with $\bI$ by fixing certain $\zs_0\in
\Z$ and let $s$ be the linear coordinate on $\bI$. Then,
$\bF_\zs(s)=\zs-s$ and, for $\zvf(s)=as+b$, we have
$$\{ as+b,\zs-s\}_\Z=-(as+b)\pa_s(\zs-s)+
(\zs-s)\pa_s(as+b)=(a\zs+b)=\zvf(\zs)$$ that proves (a). Part (b)
follows from (a) easily.
\end{pf}

Note that the vector space $\Z^\dag\simeq\wh{\Z}$ is two-dimensional but
there is no canonical basis. Instead, we have the canonical exact sequence
$$ 0\ra\bI\ra\Z^\dag\ra\bI\ra 0, $$
where the inclusion is $\bI\ni\zl\mapsto\zl 1_\Z$ and the
projection $\Z^\dag\ni\zvf\mapsto-X_\Z(\zvf)\in\bI$ gives the
`directional coefficient' of affine functions. The affine subspace
$\Z^\#$ in $\Z^\dag$ is characterized as the family of affine
functions $\zvf$ on $\Z$ for which $X_\Z(\zvf)=-1$. Similarly, the
image of $\Z$ under the isomorphism $\bF:\wh{\Z}\ra\Z^\dag$ is
characterized by $X_\Z(\zvf)=1$. The Jacobi bracket (\ref{jb})
describes the pairing between $\Z^\dag$ and $\wh{\Z}$.

\begin{theorem}\label{tp} For all $\zvf\in\Z^\dag$ and all $u\in\wh{\Z}$,
\be \{ \zvf,\bF_u\}_\Z=\lan \zvf,u\ran.\ee
\end{theorem}
\begin{pf} The theorem follows easily from (\ref{jbdu1}) by
linearity.
\end{pf}

\medskip\noindent
{\bf Remark.} We cannot add two affine scalars. However, for spaces $\Z$,
$\Z'$ of affine scalars we can introduce an equivalence relation in
$\Z\ti\Z'$ by
$$(z,z')\sim(z_1,z_1')\Leftrightarrow z-z_1=z_1'-z'$$
and interpret the equivalence class of $(z,z')$ as a sum of $z$
and $z'$. We recognize the space of such equivalence classes as
$\Z\bt\Z'$. Let us remark that this concept of addition of affine
scalars is already present in \cite{TUZ}.

\section{Affine and special affine bundles}
All above can be formulated {\it mutatis mutandis} for affine bundles
instead of affine spaces. Here $\ka=\R$ and affine bundles are smooth
bundles of affine spaces which are locally trivial, so that we pass from
one local trivialization to another using the group of affine
transformations. Since we do everything fiberwise over the same base
manifold $M$ and consider only morphisms over the identity map on the base
(if not explicitly stated otherwise), this generalization is
straightforward and we use, in principle, the same notation. For instance,
$\sV(A)$ denotes the vector bundle which is the model for an affine bundle
$\zz:A\ra M$ over a base manifold $M$. By $\Sec$ we denote the spaces of
sections, e.g. $\Sec(\zz)$ (or sometimes $\Sec(A)$) is the affine space of
sections of the affine bundle $\zz:A\ra M$. This time, however, we must
distinguish the bundles of morphisms like $\Aff_M(A_1,A_2)$,
$\Hom_M(V_1,V_2)$, etc., from their spaces of sections which consist of
particular morphisms. We will write shortly $\Aff(A_1,A_2)$ instead of
$\Sec(\Aff_M(A_1,A_2))$, etc., and $A^\dag=\Aff_M(A,\R)$ (resp.,
$V^*=\Hom_M(V,\R)$) instead of $\Aff_M(A,M\ti\R)$ (resp.,
$\Hom_M(V,M\ti\R)$) and $\Aff(A)$ (resp., $\sss{Lin}(V)$) for the space of
sections -- affine functions on $A$ (resp., linear functions on $V$).

Every section $v$ of the model vector bundle $\sV(A)$ induces a vertical
vector field $v_A$ on $A$ (called the {\it vertical lift} of $V$) being
the generator of the one parameter group of translations
$A\ni\zs_m\mapsto\zs_m+sv(m)$. Of course, $v$ is uniquely determined by
$v_A$. By a {\it special vector bundle} we understand, clearly, a vector
bundle with a distinguished nowhere vanishing section. Consequently, a
{\it special affine bundle} is an affine bundle modelled on a special
vector bundle, etc. Every special affine bundle $\bA=(A,v^0)$ carries a
distinguished vertical vector field $X_\bA=-v^0_A$, being the fundamental
vector field of the $(\R,+)$-action on $A$ induced by $v^0$, i.e. the
action. $A\ni\zs_m\mapsto\zs_m+sv^0(m)$, and thus a canonical Jacobi
structure determined by $X_\bA$. The corresponding Jacobi bracket of
smooth functions on $A$ reads
$$\{ f,g\}_\bA=fX_\bA(g)-gX_\bA(f).$$

If $V$ is a vector subbundle in the model vector bundle $\sV(A)$
of an affine bundle $A$ over $M$, then the canonical projection
$\zr:A\ra A/V$ of $A$ onto the quotient affine bundle $A/V$
defines an affine bundle structure on the total space $A$ over
$A/V$ modelled on $(A/V)\ti_M V$ (see \cite{GGU}). We will call
this affine bundle an {\it affine projection bundle} (AP-bundle)
and denote it $\AP(A,V)$.  Since $\zr$  is a morphism of affine
bundles over $M$, it makes sense to speak about the {\it affine
section bundle} $\AS(A,V)$ of $\zr$. The affine section bundle
with fibers
$$\AS(A,V)_m=\{ z_m\in\Aff(A_m/V_m;A_m):z_m\circ\zr_m=id_{A_m/V_m}\}
$$
is an affine bundle over $M$ modelled on $\Aff_M(A/V;V)$. The
space of sections of the affine bundle $\AS(A,V)$, i.e. the space
of affine sections of $\AP(A,V)$, we will denote by
$\Aff\Sec(A,V)$.

If, by chance, $A$ is a vector bundle, then we can also speak
about the {\it linear section bundle} $\LS(A,V)$ over $M$ with
fibers
$$\LS(A,V)_m=\{
u_m\in\Hom(A_m/V_m;A_m):u_m\circ\zr_m= id_{A_m/V_m}\}.
$$
This is an affine bundle over $M$ modelled on $\Hom_M(A/V;V)$. The space
of sections of the affine bundle $\LS(A,V)$, i.e. the space of linear
sections of $\AP(A,V)$, will be denoted by $\Lin\Sec(A,V)$.

Using the canonical extensions of affine maps from an affine space to linear
maps from its vector hull we get the following variant of Theorem \ref{al}.
\begin{theorem}\label{feq} The canonical embedding
$\Aff_M(A/V;A)\subset\Hom_M(\wh{A}/V;\wh{A})$ induces a canonical
identification
$$ \AS(A,V)\simeq\LS(\wh{A},V). $$
\end{theorem}
On the level of sections we denote this identification
$$\Aff\Sec(A,V)\ni\zs\mapsto\wh{\zs}\in\Lin\Sec(\wh{A},V).$$

\section{Bundles of affine values}
A particularly interesting case is that for one-dimensional
special affine bundles $\Z=(Z,v^0)$ over  $M$ which we will call
{\it bundles of affine values} (AV-bundles) and usually denote by
$\Z$. The fibers of such bundles are spaces of affine scalars
described in Section 5. The sections of an AV-bundle will play the
role of functions in our affine differential geometry that will be
developed in next sections. The model vector bundle $\sV(\Z)$ for
$\zx:\Z\ra M$ is one-dimensional and equipped with a distinguished
non-vanishing section. It is clear that this yields a canonical
identification of $\sV(\Z)$ with the trivial bundle $M\ti\R$ with
distinguished non-vanishing section represented by the constant
function $1_M$, i.e. with $M\ti\bI$. Thus the AV-bundle $\Z$
itself is {\bf trivializable}, since every section
$\zs\in\Sec(\Z)$ defines the isomorphism
$I_\zs:\Z\ra\sV(\Z)=M\ti\bI$, but {\bf not trivial}, because we
have no canonical trivialization. We insist on not introducing any
particular trivialization, since introducing it is like fixing a
frame or observer in a physical system and our approach is thought
of as a geometric framework for studying such systems in a
frame-independent way.

The sections of $\Z$ can be viewed as `functions with affine
values', since they take values in fibers of $\Z$ which are almost
reals except for the fact that we do not know where is 0, so we
can only measure the relative positions of points. The main
difference and difficulty is now that $\Sec(\Z)$ is not an algebra
nor even a vector space but only an affine space modelled on the
algebra $C^\infty(M)$ of smooth functions. In what follows, we
will identify the model bundle for an AV-bundle $\Z$ with
$M\ti\bI$. Thus we can add reals, $z_m\mapsto z_m+s$, in every
fiber $Z_m$ of $Z$, so we have a free and transitive on fibers
action of the group $(\R,+)$ on $\Z$, i.e. $\Z$ is an
$\R$-principal bundle. Let us recall that the vertical vector
field on $\Z$ which is the fundamental vector field of this action
we denote by $X_\Z$ and the corresponding vertical Jacobi bracket
on $\Z$ by $\{\cdot,\cdot\}_\Z$. The adjoint special affine bundle
$\oZ$ is represented by $Z$ with the opposite action of $\R$, i.e.
with the fundamental vector field $-X_\Z$. Conversely, it is easy
to see that every $\R$-principal bundle $Z$ carries an AV-bundle
structure. We have an obvious bundle version of Theorem \ref{tf}.

\begin{theorem}\label{ci}
There is a canonical isomorphism
\be\label{duality}\bF:\Z\ra\oZ^\#,\quad
\bF_{a_m}(a'_m)=a_m-a'_m,\ee represented by the the special affine
pairing
$$\Z\ti\oZ\ni(a_m,a'_m)\mapsto a_m-a'_m\in\bI.$$
This isomorphism extends by linearity to an isomorphism
$\bF:\wh{\Z}\ra\Z^\dag$ of special vector bundles.
\end{theorem}
$\bF:\wh{\Z}\ra\Z^\dag$ defines also a map on the level of
sections, $u\in\Sec(\wh{\Z})\mapsto\bF_u\in\Aff(\Z)$. Since
$M\ti\bI\hookrightarrow\wh{\Z}$ as $\sV(\Z)$, we can understand
$1_M$ as a section of $\wh{\Z}$ and we obtain $\bF_{1_M}=1_\Z$, so
the map $\bF$ identifies functions on $M$ with their pull-backs to
$\Z$. Moreover, for any $\zs\in\Sec(\Z)$, the function $\bF_\zs$
is an affine function on $\Z$ which is uniquely characterized by
the property that $\bF_\zs$ vanishes on the image of
$\zs\in\Sec(\Z)$ and $X_\Z(\bF_\zs)=1$. This allows us to
understand sections of $\Z$ as smooth functions $\zf$ on $\Z$ with
$X_\Z(\zf)=1$. The space of sections of $\oZ$ is identified with
the space of smooth functions on $\Z$ satisfying $X_\Z(\zf)=-1$.

An important observation is that every special affine bundle
$\bA=(A,v^0)$ gives rise to an AV-bundle. Indeed, the vector
bundle $\Aff_M(A/\la v^0\ran;\la v^0\ran)$ is  special. As the
distinguished section ${\widetilde{v}}^0$, which is constant on
fibers of $A/\la v^0\ran$ we chose
${\widetilde{v}}^0(\zr(a_m))=-v^0(m)$. Hence, $\AP(A,\lan
v^0\ran)$ is canonically an AV-bundle which we denote by
$\bAP(\bA)$. The distinguished section is chosen in such a way
that $X_{\bAP(\bA)}$ is the vertical lift $v^0_\bA$ of $v^0$, so
$\bAP(\Z)=\oZ$ and $\bAP(\Z^\#)=\Z$ for any AV-bundle $\Z$.
Moreover the map $\bF$ for $\bAP(\bA)$ is characterized by the
property that the affine function $\bF_\zs$ associated with a
section $\zs$ of $\bAP(\bA)$ satisfies $v^0_\bA(\bF_\zs)=1$ and
$\bF_\zs\circ\zs=0$. Note the isomorphism
$\ol{\bAP(\bA)}=\bAP(\ol{\bA})$. The choice of the distinguished
section in $\bAP(\bA)$ is justified by the next two theorems.

In the linear case there is an obvious identification of sections $X$ of a
vector bundle $E$ with linear functions $\zi_{E^\s}(X)$ on the dual bundle
$E^\s$, defined by the canonical pairing. If $E'$ is a submanifold of $E^\s$
(in applications $E'$ will be usually a vector or an affine subbundle), the
restriction of $\zi_{E^\s}(X)$ to $E'$ will be denoted by $\zi_{E'}(X)$. In
this notation, a section $a$ of a special affine bundle $\bA$ (regarded as a
section of $\wh{\bA}$) will give rise to a linear function
$\zi_{\bA^\dag}(a)$ on $\bA^\dag$ and an affine function $\zi_{\bA^\#}(a)$ on
the affine subbundle $\bA^\#$ of $\bA^\dag$. Denote the map $\bF$ for the
AV-bundle $\bAP(\bA^\dag)$ (resp., $\bAP(\bA^\#)$) by $\bF^\dag$ (resp.,
$\bF^\#$).

For a special affine bundle (resp., a special vector bundle) $\bA=(A,v^0)$
denote $\AS(A,\langle v^0\rangle)$ by $\bAS(\bA)$ (resp., $\LS(A,\lan v^0\ran)$
by $\bLS(\bA)$). The spaces of sections of these bundles we denote simply
$\Aff\Sec(\bA)$ and $\Lin\Sec(\bA)$, respectively. Since the section
${\widetilde{v}}^0$ is affine, also $\bAS(\bA)$ is canonically a special affine
bundle. In the case when $\bA$ is a special vector bundle the affine bundle
$\bLS(\bA)$ is not canonically special, since the section ${\wt{v}}^0$ is not
linear. However, in the case when $\bA$ is a bispecial vector bundle with the
distinguished section $\zf^0$ of $\bA^\s$, $\la\zf^0,v^0\ran=0$, then also
$\bLS(\bA)$ is special affine with the distinguished section
$\wt{v_{\zf^0}}\in\Hom(A/\la v^0\ran;\la v^0\ran)$,
$$\wt{v_{\zf^0}}(\zr(a_m))=-\la\zf^0(m),a_m\ran v^0(m).$$

\begin{theorem}\label{1} \cite{GGU} There is a canonical isomorphism
of affine bundles
$$A\simeq{\bLS}(A^\dag),\quad a_m\mapsto\wh{\zs}_{a_m},$$
where
$$\wh{\zs}_{a_m}([\zf_m])=\zf_m-\zf_m(a_m)1_A(m).$$
In other words, for any section $a$ of $A$,
$$
\bF^\dag_{\wh{\zs}_a}=\zi_{A^\dag}(a).
$$
The corresponding isomorphism of the model vector bundles takes the form
$$\sV(A)_m\ni X_m\ \longleftrightarrow\ -\zi^\dag_{X_m}\in(A^\dag/\la
1_A\ran)^\s,$$ where $\zi^\dag_{X_m}([\zf_m])=(\zf_m)_\sv(X_m)$.
\end{theorem}
Note that the above theorem is an affine version of the well-known
fact that sections of a vector bundle $E$ over $M$ can be
identified with linear (along fibers) functions on the dual
$E^\s$, i.e. with linear sections of the bundle $E^\s\ti\R$ over
$E^\s$. We can extend this identification to special affine
bundles as follows.
\begin{theorem}\label{oto} For a special affine bundle $\bA=(A,v^0)$ there is
a canonical identification of special affine bundles
\bea&{\bA}\simeq{\bAS}(\bA^\#)\simeq{\bLS}({\bA}^\dag)\cr
&a_m\longleftrightarrow{\zs}_{a_m}\longleftrightarrow \wh{\zs}_{a_m}.\eea

On the level of sections it takes the form
\bea&\Sec({\bA}) \simeq {\Aff\Sec}(\bA^\#)\simeq {\Lin\Sec}({\bA}^\dag),\cr
&a\ \ \longleftrightarrow {\zs}_a\longleftrightarrow\ \wh{\zs}_a,
\eea
where $$ \bF^\dag_{\wh{\zs}_a}=\zi_{\bA^\dag}(a),
$$
$$ \bF^\#_{{\zs}_a}=\zi_{\bA^\#}(a).
$$
This identification leads to the obvious identification of the
corresponding model vector bundles
$$\sV(\bA)=\Aff_M(\bA^\#/\lan 1_\bA\ran;\R)=\Hom_M(\bA^\dag/\lan 1_\bA\ran;\R)
(=(\bA^\dag/\lan 1_\bA\ran)^\s)$$ taking on sections the form
$$X\longleftrightarrow-\zi_X^\#\longleftrightarrow-\zi_X^\dag,$$
where linear functions $\zi_X^\dag$ and affine functions $\zi_X^\#$ on
$\bA^\dag/\langle 1_\bA\rangle$ and $\bA^\#/\langle 1_\bA\rangle$,
respectively, are the projections of linear functions $\zi_{\bA^\dag}(X)$
on $\bA^\dag$ and  affine functions $\zi_{\bA^\#}(X)$ on $\bA^\#$,
respectively.
\end{theorem}
\begin{pf} The proof that these bundles are canonically isomorphic is just
the combination of Theorems \ref{1} and \ref{feq}. That the distinguished
sections are preserved follows from
$$
\bF^\dag_{\wh{\zs}_{a+v^0}}=\zi_{\bA^\dag}(a+v^0)=\bF^\dag_{\wh{\zs}_a}+\zi_{\bA^\dag}(v^0),
$$
$$ \bF^\#_{{\zs}_{a+v^0}}=\zi_{\bA^\#}(a+v^0)=\bF^\#_{{\zs}_a}+1.
$$
\end{pf}

\begin{corollary}\label{L} For an affine bundle $A$ and an AV-bundle $\Z$
over $M$ there is a canonical identification
$$\Aff_M(A;\Z)\simeq
\bA^\dag\bt_M\Z.$$
\end{corollary}
\begin{pf} Observe first that $A\ati_M\Z$ is canonically a special
affine bundle and the identification {\it mapping $\leftrightarrow$ graph}
induces the identification of $\Aff_M(A;\Z)$ with the space
${\bAS}(A\ati_M\oZ)$ of affine sections of the associated AP-bundle
$A\ati_M\oZ$ over $A$ . The latter is, due to the above theorem, canonically
identified with the special affine bundle $(A\ati_M\ol{\Z})^\#$ over $M$. In
view of Theorem \ref{bt} and Theorem \ref{ci}
$$
(A\ati_M\ol{\Z})^\#\simeq \bA^\dag\bt_M\ol{\Z}^\#\simeq
\bA^\dag\bt_M\Z.$$
\end{pf}

We will end up this section with presenting the above concepts in
local coordinates. First of all, for a special vector bundle
$\sv(\zh)\colon\bV=(V,v^0)\ra M$ we choose a coordinate
neighborhood $U$ in $M$ with coordinates $x=(x^b)$ and a basis
$(v^1,\dots,v^k,v^0)$ of local sections over $U$ which contains
the distinguished $v^0$. On fibers over $U$ we have then the
associated linear coordinates $(y,s)=(y_1,\dots,y_k,s)$, so the
coordinates $(x,y,s)$ on $(\sv(\zh))^{-1}(U)$. We will call such
local coordinates on $\bV$ {\it linear coordinates}. These
coordinates can serve as coordinates on $\zh^{-1}(U)$ for any
special affine bundle $\zh:\bA=(A,v^0)\ra M$ modelled on
$\sv(\zh)$ if we use the isomorphism of affine bundles $I_\zs:A\ra
V=\sV(A)$ determined by a section $\zs\in\Sec(A)$. Such
coordinates will be called {\it local affine coordinates on
$\bA$}. The change of the section $\zs$ results in the
transformation of coordinates by a translation
$(x,y,s)\mapsto(x,y+f(x),s+g(x))$, so that objects of affine
differential geometry should be defined in local coordinates
invariantly with respect to this change of coordinates. On
$\wh{\bA}$ we have linear coordinates $(x,y,z,s)$ such that $A$ is
characterized by the equation $z=1$. The canonical vector field on
$\bA$ has the expression $X_\bA=-\pa_s$. Affine functions on $A$
have the form
$$\zf(x,y,s)=\za^i(x)y_i+\zg(x)s+\zb(x)$$
and correspond to linear functions
$$\wh{\zf}(x,y,z,s)=\za^i(x)y_i+\zg(x)s+\zb(x)z$$
on $\wh{\bA}$. Hence, $(x^b,\za^1,\dots,\za^k,\zb,\zg)$ represent
coordinates on $\bA^\dag$. The distinguished section is
$1_\bA(x)=(x,0,1,0)$. The affine subspace $\bA^\#$ in $\bA^\dag$ is
characterized by $\zg=1$.

If $\bA=\Z$ is an AV-bundle then the coordinates $y$ are lacking and the
affine function corresponding to the section $\zs:s=\zs(x)$ is
$\bF_\zs(x,s)=\zs(x)-s$. For the particular case of the AV-bundles
$\bAP(\bA^\#)$ and $\bAP({\bA}^\dag)$ induced by a special affine bundle
$\bA$ we have coordinate expressions $(x,\za,\zb)\mapsto(x,\za)$ and
$(x,\za,\zb,\zg)\mapsto(x,\za,\zg)$, respectively. The distinguished
sections are described by the equation $\zb=-1$. The canonical pairing
between $\wh{\bA}$ and $\bA^\dag$ is
$$\lan(x,y,z,s),(x,\za,\zb,\zg)\ran=y_i\za^i+z\zb+s\zg,$$
so that the canonical pairing between $\bA$ and $\bA^\#$ reads
$$\lan(x,y,s),(x,\za,\zb)\ran_{sa}=
\lan(x,y,1,s),(x,\za,\zb,1)\ran=y_i\za^i+\zb+s.$$ In other words,
$\zi_{\bA^\#}(a)(x,\za,\zb)=y_i(x)\za^i+\zb+s(x)$ for a section
$a(x)=(x,y(x),s(x))$ of $\bA$.

Affine (resp., linear) sections of the bundles $\bAP(\bA^\#)$ and
$\bAP({\bA}^\dag)$ have the form
$$\zs(x,\za)=(x,\za,y_i(x)\za^i+s(x))\quad \textrm{and}\quad
\wh{\zs}(x,\za,\zg)=(x,\za,\zg,y_i(x)\za^i(x)+s(x)\zg),
$$
respectively. The associated affine function
$\bF_\zs^\#=i_{\bA^\#}(\ol{a})$ on $\bA^\#$ reads
$$\bF_\zs(x,\za,\zb)=\zb-y_i(x)\za^i-s(x)=\lan(x,\za,\zb),
(x,-y(x),-s(x))\ran_{sa}
$$
and corresponds to the section $\ol{a}(x)=(x,-y(x),-s(x))$ of
${\bA}$. Conversely, the section $a(x)$ corresponds, by
definition, to the affine section
$$\zs_a(x,\za)=(x,\za,\zb-\lan(x,y(x),s(x)),(x,\za,\zb)\ran_{sa})=
(x,\za,-y_i(x)\za^i-s(x))$$ of $\bA^\#$.

\section{AV-differential geometry: the phase and the contact bundles}
The standard Cartan calculus of differential forms is based on the
algebra of differentiable functions on a manifold $M$. We will
start to build {\it AV-differential geometry} where for the Cartan
calculus we replace functions by sections of an AV-bundle $\zz
\colon \Z\rightarrow M$ modelled on the trivial bundle
$pr_M:M\times\bI$. This is our starting object whose sections
$\Sec(\Z)$ replace the sections of $M\ti\R$, i.e. smooth functions
$C^\infty(M)$ on $M$ in the standard differential geometry. This
chapter is based on \cite{TU,Ur2}, where AV-analogs of the
cotangent and contact bundles $\sT^\s M$ and $\sT^\s M\ti\R$ have
been introduced.

One builds an AV-analog of the cotangent bundle $\sT^\s M$ as
follows. Let us define an equivalence relation in the set of all
pairs $(m,\zs )$, where $m$ is a point in $M$ and $\sigma $ is a
section of $\zz $. Two pairs $(m,\zs )$ and $(m',\zs ')$ are
equivalent if $m' = m$ and $\xd(\zs ' - \zs)(m) = 0$.  We have
identified the section $\zs ' - \zs $ of $pr_M$ with a function on
$M$ for the purpose of evaluating the differential $\xd(\zs ' -
\zs)(m)$.   We denote by $\sP\Z$ the set of equivalence classes.
The class of $(m,\zs )$ will be denoted by $\bd\sigma (m)$ or by
$\bd_m\zs$ and will be called the {\it differential} of $\zs $ at
$m$. We will write $\bd$ for the affine exterior differential to
distinguish it from the standard  $\xd$. We define a mapping
${\sP}\zeta \colon \sP \Z \rightarrow M$ by $\sP\zz (\bd\zs (m)) =
m$. The bundle $\sP\zz$ is canonically an affine bundle modelled
on $\zp _M \colon \sT^{\textstyle *} M \rightarrow M$ with the
affine structure
$$\bd\zs _2(m)- \bd\zs _1(m) =
\xd(\zs _2 - \zs _1)(m).
$$
This affine bundle is called the {\it phase bundle} of $\zz$. A section of
$\sP\zz$ will be called an {\it affine 1-form}.

Let $\za\colon M\rightarrow \sP\Z$ be an affine 1-form and let $\zs$ be a
section of $\zz$. The differential $\xd_m(\za -\bd \zs)$ does not depend on
the choice of $\zs$ and will be called the {\it differential of $\za$ at
$m$}. We will denote it by $\bd\za(m)$ or by $\bd_m\za$. The differential of
an affine 1-form $\za\in\Sec(\sP\Z)$ is an ordinary 2-form
$\bd\za\in\zW^2(M)$. The corresponding {\it affine de Rham complex} looks now
like \be\label{DR}\Sec(\Z)\bdra\Sec(\sP\Z)\bdra\zW^2(M)\dra\zW^3(M)\dra\dots
\ee and consists of affine maps. This is an {\it affine complex} in this
sense that its linear part is a complex of linear maps, so that we
can define the corresponding cohomology. The linear part of
(\ref{DR}) is a part of the standard de Rham complex (without its
beginning consisting of the inclusion of $\R$ into $C^\infty(M)$).
However, note that the cohomology of (\ref{DR}) can be defined
without refereing to its linear part. Indeed, the problem is only
with the first and the second cohomology space, since the rest is
the standard de Rham complex. Denote the kernel (the inverse image
of $\{ 0\}$) and the image of the affine map
$\bd:\Sec(\sP\Z)\ra\zW^2(M)$ by $Z_1$ and $B_2$, respectively.
$Z_1$ is an affine subspace of $\Sec(\sP\Z)$ and $B_2$ is a vector
subspace of the kernel $Z_2$ of $\xd:\zW^2(M)\ra\zW^3(M)$.
Moreover, the image $B_1$ of $\bd:\Sec(\Z)\ra\Sec(\sP\Z)$ is an
affine subspace in $Z_1$. But the quotients of affine spaces are
vector spaces, so that $H^1=Z_1/B_1$ and $H^2=Z_2/B_2$ are vector
spaces. It is easy to see that we got nothing but the lacking
first and second de Rham cohomology.

Recall that the bundle $\Z$ can be considered as a principal bundle with
the structure group $(\R,+)$ and the fundamental vector field of this
action we have denoted by $X_\Z$. Let us observe now that $\sP\Z$
represents the principal connection bundle of $\Z$, i.e. there is a
canonical identification of the affine space $\Sec(\sP\Z)$ of sections of
$\sP\Z$ with the affine space $\sss{PConn}(\Z)$ of principal connections
in $\Z$. We will identify the space of principal connections with the
space of connection 1-forms. In other words, $\sss{PConn}(\Z)$ consists of
those 1-forms $\zn$ on $\Z$ which are $X_\Z$-invariant and $\zn(X_\Z)=1$.
Since we can add to $\zn$ the pull-backs of 1-forms on $M$, both affine
spaces are modelled on the space $\zW^1(M)$ of 1-forms on $M$.

\begin{theorem} There is a canonical isomorphism of affine spaces
$\sF:\Sec(\sP\Z)\ra\sss{PConn}(\Z)$, with linear part being the identity
on $\zW^1(M)$, such that for any section $\zs$ of $\Z$
$$\sF_{\bd\zs}=\xd
\bF_\zs,
$$
where $\bF_\zs(\zs')=\zs-\zs'$. Moreover,
$$\xd\sF_\za=\zz^\s(\bd\za),$$
so that the 2-form $\bd\za\in\zW^2(M)$ is the curvature form of the
connection $\za\in\Sec(\sP\Z)$.
\end{theorem}
\begin{pf} Indeed, since $X_\Z(\bF_\zs)=1$, $\xd
\bF_{\zs}\in\sss{PConn}(\Z)$. Moreover, for $f\in C^\infty(M)$,
$$\sF_{\bd(\zs+f)}=\xd(\bF_\zs+f\circ\zz),$$ so that
$$\sF_{\bd\zs+\xd f}=\sF_{\bd\zs}+\zz^*\xd f$$
and we can define $\sF_\za$ for arbitrary $\za\in\Sec(\sP\Z)$ by
$\sF_\za=\sF_{\bd\zs}+\zz^\s(\za-\bd\zs)$. Finally,
$$\xd\sF_\za=\xd(\sF_{\bd\zs}+\zz^\s(\za-\bd\zs))=\zz^\s\xd(\za-\bd\zs)=
\zz^\s\bd\za.$$ Conversely, if $\zn\in\sss{PConn}(\Z)$, then $\zn-\bd\zs$
is a vertical and $X_\Z$-invariant one-form on $\Z$ for any section $\zs$
of $\Z$, thus $\zn-\bd\zs=\zz^\s(\zm)$ for certain one-form $\zm$ on $M$.
Then, $\zn=\sF_\za$ for $\za=\bd\zs+\zm$.
\end{pf}

\medskip\noindent
In local affine coordinates $(x^a,s)$ on $\Z$ we have
$$\sF_{\za_a(x)\xd x^a}=\za_a(x)\xd x^a-\xd s.$$

\medskip
As we have noticed, there is a distinguished affine space $Z_1$ of
{\it closed} affine 1-forms. It turns out that, like in the case
of the cotangent bundle, they can be defined intrinsically as
those sections of $\sP\Z$ whose images are lagrangian submanifolds
with respect to a {\it canonical symplectic structure} on $\sP\Z$
which is defined as follows.

For a chosen section $\zs$ of $\zz$ we have isomorphisms
      \bea\nn
   &I_\zs\colon\Z\rightarrow M\times \R,  \\
      &I_{\bd\zs}\colon \sP\Z\rightarrow \sT^{\textstyle *} M,
                                       \label{Fa29}\eea
   and for two sections $\zs, \zs'$ the mappings $I_{\bd\zs}$ and
   $I_{\bd\zs'}$ differ by translation by $\xd(\zs -\zs')$, i.e.
      \bea\nn
   &I_{\bd\zs'}\circ I_{\bd\zs}^{-1}\colon \sT^{\textstyle *} M\rightarrow
\sT^{\textstyle *} M  \\
   &\colon\za_m \mapsto \za_m +\xd(\zs -\zs')(m).
                                       \label{Fa30}\eea
Now we use an affine property of the canonical symplectic form
$\zw_M$ on the cotangent bundle: {\bf translations in $\sT^\s M$
by a closed 1-form are symplectomorphisms}, to conclude that the
two-form $I_{\bd\zs}^{\textstyle *} \zw_M$, where $\zw_M$ is the
canonical symplectic form on $\sT^{\textstyle *} M$, does not
depend on the choice of $\zs$ and therefore it is a canonical
symplectic form on $\sP\Z$. We will denote this form by $\zw_\bZ$.

\begin{theorem} An affine 1-form $\za\in\Sec(\sP\Z)$ is closed if
and only if $\za(M)$ is a lagrangian submanifold of
$(\sP\Z,\zw_\Z)$.
\end{theorem}
\begin{pf} Consider a section $\zs\in\Sec(\Z)$ and the
corresponding isomorphism of affine bundles $I_{\bd\zs}\colon
\sP\Z\rightarrow \sT^{\textstyle *} M$. With respect to this
isomorphism any affine 1-form $\za$ corresponds to the true 1-form
$\za-\bd\zs$ on $M$: $I_{\bd\zs}(\za(m))=\za(m)-\bd\zs(m)$.
According to the well-known characterization, $\za-\bd\zs$ is
closed if and only if $(\za-\bd\zs)(M)$ is a lagrangian submanifod
in $(\sT^\s M,\zw_M)$ so if and only if
$$\za(M)=I_{\bd\zs}^{-1}((\za-\bd\zs)(M))$$  is a lagrangian
submanifold of $(\sP\Z,\zw_\Z)$, since $I_{\bd\zs}$ is a
symplectomorphism. But, by definition, $\xd(\za-\bd\zs)=0$ if and
only if $\bd\za=0$.
\end{pf}

\medskip\noindent
{\bf Remark.} It is obvious that $\sP\Z$ and $\sP \overline{\Z}$
are equal as manifolds. Let $\zs$ be a section of $\Z$. The same
mapping interpreted as a section of  $\overline{\Z}$ will be
denoted by $\overline{\zs}$. Since $\zs -\zs' = \overline{\zs}' -
\overline{\zs}$, the isomorphisms $I_{\bd\zs}\colon \sP \Z
\rightarrow \sT^{\textstyle *} M$ and $I_{\bd\overline{\zs}}\colon
\sP \overline{\Z}\rightarrow \sT^{\textstyle *} M $ are related by
$I_{\bd\zs} = - I_{\bd\overline{\zs}}$. It follows that
   $$\zw_{\overline{\bZ}} = I_{\bd\overline{\zs}}^{\textstyle *} \zw_M = -
I_{\bd\zs}^{\textstyle *} \zw_M =-\zw_{\bZ}.$$

\medskip\noindent
There is no canonical Liouville 1-form on $\sP\Z$ (in the standard
sense) which is the potential of the canonical symplectic form
$\zw_\Z$ but there is such a form in the affine sense. To define
this Liouville 1-form let us build another canonical affine bundle
out of $\Z$.

   We define another equivalence relation in the set of all pairs $(m,\zs
)$. Two pairs $(m,\zs )$ and $(m',\zs ')$ are equivalent if $m' =
m$, $\zs(m) = \zs'(m)$, and $\xd(\zs ' - \zs)(m) = 0$.  We can
identify the equivalence class of $(m,\zs)$ with the first jet of
the section $\zs$ with the source point $m$.  We denote by $\sC\Z$
the set of equivalence classes. The class of $(m,\zs )$ will be
denoted by $\xc\sigma (m)$ or by $\xc_m\zs$ and will be called the
{\it contact element} of $\zs $ at $m$. We define a fiber bundle
structure over $M$ defining the projection ${\sC}\zeta \colon
\sC\Z \rightarrow M$ by $\sC\zz (\xc\zs (m)) = m$. In other words,
$\sC\Z$ is the first-jet bundle $j^1(\zz)$ of $\zz$. This fiber
bundle is canonically an affine bundle modelled on $\zg _M \colon
\sT^{\textstyle *} M \op \R \rightarrow M$ with the affine
structure defined by
  $$\xc\zs _2(m)- \xc\zs _1(m) = (\xd(\zs _2 - \zs _1)(m),
\zs_2(m)- \zs_1(m)).$$ This affine bundle is called the {\it
contact bundle} of $\zz$. The pair $(\sC\Z,(0, 1_M))$ is a special
affine bundle. It is easy to see $\sC\Z=\sP\Z\ati_M\Z$ and that
$\sC\Z/\lan(0, 1_M)\ran$ is canonically isomorphic to $\sP\Z$ (we
just identify the points $m$ in the equivalence relation) so we
have the associated AV-bundle with the canonical projection
$\zz_{\sC \bZ}:\sC\Z\ra\sP\Z$.

There is also  a canonical projection
   \be\zm\colon \sC\Z \rightarrow\Z,\quad \zm(\xc\zs(m))=\zs(m),
                              \label{Fa32}\ee
   which is a morphism of special affine bundles $\zz_{\sC \bZ}\colon  \sC\Z
\rightarrow \sP\Z$ and $\zz \colon\Z\rightarrow M$ over the
projection $\zz:\sP\Z\ra M$ on the level of base manifolds and
there is a well-defined pull-back of sections of $\zz$ to sections
of $\zz_{\sC\bZ}$. Now we can define a section $\zvy_\bZ$ of
$\sP\zz_{\sC\bZ}\colon \sP \sC\Z \rightarrow \sP\Z$ by
      \be \zvy_\bZ(p) = \bd_m \zm^{\textstyle *} \zs,
                                                       \label{Fa33}\ee
$\sP\zz(p)=m$, where $\zs$ is a section of $\zz$ which represents
$p\in \sP_m\Z$. In other words, for any section $\zs\in\Sec(\Z)$
$$\zvy_\Z(\bd\zs(m))=\bd(\zm^\s\zs)(\bd\zs(m)).$$
The affine 1-form $\zvy_\bZ$ is called the {\it Liouville affine form} of
$\sC\Z$ and defines the {\it canonical contact structure} of $\sC\Z$. This
affine 1-form over $\sP\Z$ is a potential for the canonical symplectic
form on $\sP\Z$.

\begin{theorem} $\bd \zvy_\bZ=\zw_\Z$.
\end{theorem}
\begin{pf} Let us take a section $\zs_0$ of $\Z$. Using the
identification $I_{\xc\zs_0}:\sC\Z\ra\sT^\s M\ti\R$ we identify
the AV-bundle $\zz_{\sC\Z}:\sC\Z\ra\sP\Z$ with the trivial bundle
$pr_{\sT^\s M}:\sT^\s M\ti\R\ra\sT^\s M$. With this identification
sections of $\zz_{\sC\Z}$ are functions on $\sT^\s M$ and
$\zm^\s\zs_0=0$ is a distinguished section of $\zz_{\sC\Z}$, so
that sections of $\sP\sC\Z$ are standard 1-forms on $\sT^\s M$
with the standard de Rham differential. Moreover, sections of $\Z$
are represented by functions on $M$, $\zm^\s$ is represented by
$\zp_M^\s$, and the symplectic form $\zw_\Z$ is represented by the
standard symplectic form $\zw_M$ on $\sT^\s M$.

Take a section $\zs$ of $\Z$ understood as a function on $M$. By
definition, $\zvy_\Z(\xd\zs(m))=\xd(\zp_M^\s\zs)(\xd\zs(m))$ which
means that $\zvy_\Z$ is represented by the true Liouville 1-form
$\zvy_M$ on $\sT^\s M$. Hence, $\bd\zvy_\Z=\zw_\Z$.
\end{pf}

\medskip\noindent
{\bf Remark.} Note that the above proof does not imply that we can
define a true canonical Liouville 1-form on $\sP\Z$. Indeed, it is
easy to see that the change of the initial section  $\zs_0$ into
$\zs_0'$ with $\zs'_0=\zs_0+f$ results, for the trivialization
given by $\zs_0$, in translation of the Liouville 1-form:
$\zvy_M\mapsto \zvy_M-\zp_M^\s\xd f$. Thus, the true Liouville
1-form on $\sT^\s M$ has no affine meaning (but its exterior
derivative has such a meaning), since it is not invariant by
translations by $\xd\zp_M^\s(f)$). We put a geometrical meaning to
the transformation rules of the Liouville 1-form defining its
affine version. This explains perhaps better what an affine 1-form
is.

\medskip
The affine Liouville 1-form $\zvy_\Z$ can be interpreted as a
canonical principal connection on the principal bundle
$\zz_{\sC\Z}:\sC\Z\ra\sP\Z$, thus as a canonical 1-form
$\zh_{\sC\Z}=\sF_{\zvy_{\Z}}$ on $\sC\Z$. In any trivialization
$I_{\xc\zs}:\sC\Z\ra\sT^\s M\ti\R$ and the standard Darboux
coordinates $(x^a,p_b,s)$ on $\sT^\s M\ti\R$ the affine Liouville
1-form has the standard expression $\zvy_\Z=p_a\xd x^a$, so
$\zh_{\sC\Z}$ looks like the canonical contact form:
$\zh_{\sC\Z}=p_a\xd x^a-\xd s$. It can be also seen directly that
the canonical contact form $\zh_{M}=p_a\xd x^a-\xd s$ on $\sT^\s
M\ti\R$ is affine in this sense that it is invariant with respect
to translations of $\sT^\s M\ti\R$ by first jets of functions . We
will call $\zh_{\sC\Z}$ the {\it canonical contact form} on
$\sC\Z$. Like every contact form, it induces a (contact) Jacobi
bracket $\{\cdot,\cdot\}_{\sC\Z}$ on $\sC\Z$ which in the above
local coordinates reads

\be\label{cjb}\{ f,g\}_{\sC\Z}=\frac{\pa f}{\pa{p_a}}\frac{\pa g}{\pa{x^a}}-
\frac{\pa g}{\pa{p_a}}\frac{\pa f}{\pa{x^a}} + (p_a\frac{\pa
f}{\pa{p_a}}-f)\frac{\pa g}{\pa s}- (p_a\frac{\pa
g}{\pa{p_a}}-g)\frac{\pa f}{\pa s}.\ee The corresponding Jacobi
structure on $\sC\Z$ in this trivialization takes the form
$(\zL_M+\zD_{\sT^\s M}\we\pa_s,-\pa_s)$, where
$\zL_M=\pa_{p_a}\we\pa_{x^a}$ is the canonical Poisson tensor on
$\sT^\s M$ associated with the syplectic form $\zw_M$ and
$\zD_{\sT^\s M}=p_a\pa_{p_a}$ is the Liouville vector field on
$\sT^\s M$ (both regarded as tensor fields on $\sT^\s M\ti\R$).
Again, this Jacobi structure has an affine flavor, since the
tensors are invariant with respect to translations in $\sT^\s
M\ti\R$ by first jets of functions.

\section{Lie algebroids associated with AV-bundles}

The principal bundle structure of $\Z$ represented by the vector
field $X_\Z$ induces additional structures on functions, vector
fields and, in general, differential operators on $Z$. For any
manifold $N$ denote by $\X(N)$ (resp., $\de(N)$) the space of all
vector fields (resp., the space of all linear first-order
differential operators) on $N$, i.e. acting on $C^\infty(N)$.
Clearly, $\X(N)=\Sec(\sT N)$ and $\de(N)$ is the space of sections
of the bundle $\sL N=\sT N \oplus\R$ of linear first-order
differential operators on $N$ with the obvious action of
$(X+h)\in\Sec(\sL N)=\X(N)\op C^\infty(N)$ on functions on $N$
given by $(X+h)(f)=X(f)+hf$.

Let us fix an AV-bundle $\Z$ over $M$. The space $Pol^n(\Z)$ of
polynomials of order $\le n$ on $\Z$ is defined as the space of
those smooth functions $f$ on $\Z$ for which $X_\Z^{n+1}(f)=0$, so
that we have the filtration $Pol(\Z)=\bigcup_nPol^n(\Z)$ of the
algebra of all polynomials. Note that in the affine case we have
only the filtration and no canonical graduation of $Pol(\Z)$. In
particular, the space $Pol^0(\Z)$ is just the algebra $Bas(\Z)$ of
basic functions on $\Z$, i.e. functions that are constant along
fibers (it will be often identified with the algebra of smooth
functions on $M$) and $Pol^1(\Z)$ is the space $\Aff(\Z)$ of
affine (along fibers) functions on $\Z$.

We have also natural subalgebras of the Lie algebra $\X(\Z)$ of
all vector fields on $\Z$. The Lie algebra $\wt{\X}(\Z)$ of
invariant vector fields on $\Z$ consists of those vector fields
$X$ for which $[X_\Z,X]=0$. It is easy to see that, in local
affine coordinates $(x^a,s)$ on $\Z$, invariant vector fields have
the form
$$ X=f_a(x)\pa_{x^a}-g(x)\pa_s, $$
where the functions $f_a,g$ are basic. Vector fields from $\wt{\X}(\Z)$
can be viewed as sections of the vector bundle
$\widetilde{\sT}\Z=\sT\Z/\R$ which is the vector bundle over $M$ of orbits
of the tangent lift $\zvf_\s$ of the $(\R,+)$ action $\zvf$ on $\Z$. Since
the vector field $X_\Z$ is invariant, it can be understood as a
distinguished section of $\widetilde{\sT}\Z$, so $\widetilde{\sT}\Z$ is
canonically a special vector bundle. This is just the Atiyah vector bundle
(and canonically a Lie algebroid) associated with the $\R$-principal
bundle $\Z$.

There is another natural subalgebra of the Lie algebra $\X(\Z)$ of
all smooth vector fields on $\Z$, namely the subalgebra
$\X_{ah}(\Z)$ of {\it affine-homogeneous vector fields}, i.e.
those vector fields $X$ which preserve the filtration:
$X(Pol^n(\Z))\subset Pol^n(\Z)$. Of course,
$\wt{\X}(\Z)\subset\X_{ah}(\Z)$, since invariant vector fields
lower the filtration: $X(Pol^n(\Z))\subset Pol^{n-1}(\Z)$. Again,
the space $\X_{ah}(\Z)$ is the space of sections of certain vector
bundle $\breve{\sL}\Z$ over $M$ which can be identified with the
bundle $\wt{\sT}\Z\svop_M\Z^\dag=(\wt{\sT}\Z\op_M\Z^\dag)/\lan
X_\Z-1_\Z\ran$, i.e. the special direct sum of the special vector
bundles $\wt{\sT}\Z$ and $\Z^\dag$. Indeed, it is easy to see that
the class $Y\svop\zf\in\Sec(\wt{\sT}\Z\svop_M\Z^\dag)$, where
$Y=f_a(x)\pa_{x^a}-g(x)\pa_s$ and $\zf(x,s)=\za(x)s+\zb_0(x)$, is
represented by $f_a(x)\pa_{x^a}-(\za(x)s+\zb(x))\pa_s$ with
$\zb(x)=\zb_0(x)+g(x)$, so the vector field
${\obD}_{X\svop\zf}=X+\zf X_\Z$ gives such an identification.

Similarly, there is a natural subalgebra of the Lie algebra
$\de(\Z)=\de(Z)$, the subalgebra $\de_{ah}(\Z)$ of {\it
affine-homogeneous first-order differential operators}, consisting
of those $D\in\de(\Z)$ which preserve the filtration:
$D(Pol^n(\Z))\subset Pol^n(\Z)$. Note that $\de_{ah}(\Z)$ is
canonically a $Bas(\Z)\simeq C^\infty(M)$-module.

It is easy to see the following.

\begin{proposition} There is a canonical splitting $\de_{ah}(\Z)=
\X_{ah}(\Z)\op Bas(\Z)$. Moreover, a vector field $X$ on $\Z$ is
affine-homogeneous if and only if $[X_\Z,[X_\Z,X]]=0$ and
$[X_\Z,X]$ is vertical. In local affine coordinates $(x^a,s)$ on
$\Z$, affine-homogeneous vector fields have precisely the form
$$
X=f_a(x)\pa_{x^a}-(\za(x)s+\zb(x))\pa_s,
$$
and affine-homogeneous first-order differential operators the form
$$
D=f_a(x)\pa_{x^a}-(\za(x)s+\zb(x))\pa_s+\zg(x).
$$
\end{proposition}
Note that the vector fields from $\X_{ah}(\Z)$ are projectable and
the vector field $X=f_a(x)\pa_{x^a}-(\za(x)s+\zb(x))\pa_s$
projects onto the vector field $\oX=f_a(x)\pa_{x^a}$ on $M$.
Before finding an appropriate bundle whose sections form
$\de_{ah}(\Z)$ let us observe that the canonical Jacobi bracket
$\{\cdot,\cdot\}_\Z$ applied to affine functions
$\zf,\zc\in\Aff(\Z)$ gives a basic function. Indeed, since
$X_\Z^2(\zf)=X_\Z^2(\zc)=0$, we have
$$X_\Z(\{\zf,\zc\}_\Z)=X_\Z(\zf X_\Z(\zc)-\zc X_\Z(\zf))=0.$$
Recall that the map $\bF$ identifies sections of $\wh{\Z}$ with
$\Aff(\Z)$. In particular, we can identify sections $\zs$ of $\Z$ with
affine functions $\bF_\zs$ which satisfy $X_\Z(\bF_\zs)=1$, so that
$\Sec(\Z)\simeq\{\zf\in\Aff(\Z):X_\Z(\zf)=1_\Z\}$. We have the following
bundle version of Theorem \ref{tp} obtained just fiberwise.
\begin{theorem}\label{tp1} For all $\zvf\in\Aff(\Z)=\Sec(\Z^\dag)$ and all
$u\in\Sec(\wh{\Z})=\Sec((\Z^\dag)^\s)$,
\be \{\zvf,\bF_u\}_\Z=\lan \zvf,u\ran.\ee
\end{theorem}
There is a `hamiltonian map'
$$\Aff(\Z)\ni\zf\mapsto\bD_\zf=\zf
X_\Z-X_\Z(\zf)=\{\zf,\cdot\}_\Z\in\de(\Z)
$$
with the property
$$\bD_\zf(\bF_u)=\{\zf,\bF_u\}_\Z=\lan\zf,u\ran
$$
for $\zf\in\Aff(\Z)$, $u\in\Sec(\wh{\Z})$. Therefore we can
consider $\Z^\dag$ as embedded in $\de_{ah}(\Z)$. For a section
$\zs$ of $\wh{\Z}$ we will write shortly $\bD_\zs$ instead of
$\bD_{\bF_\zs}$.

In local affine coordinates, the Jacobi bracket of affine functions on
$\Z$ takes the form
$$\{\za(x)s+\zb(x),\za'(x)s+\zb'(x)\}_\Z=\za(x)\zb'(x)-\za'(x)\zb(x)
$$
and the differential operator associated to
$\za(x)s+\zb(x)\in\Aff(\Z)$ reads
\be\label{D}
\bD_{\za(x)s+\zb(x)}=\za(x)-(\za(x)s+\zb(x))\pa_s.
\ee
Now, we can extend the map $\bD$ to sections of the bundle
$\sR\Z=\breve{\sL}\Z\svop_M\Z^\dag=\wt{\sT}\Z\svop_M\Z^\dag\svop_M\Z^\dag$
by
$$ \bD_{X\svop\zf\svop\zc}=X+\zf X_\Z+\bD_\zc.
$$
It is easy to see that this gives the identification of sections of $\sR\Z$
with $\de_{ah}(\Z)$. In local affine coordinates,
$$\bD_R=
f_a(x)\pa_{x^a}-((\za(x)+\za'(x))s+g(x)+\zb(x)+\zb'(x))\pa_s+\za(x),
$$
where $R={(f_a(x)\pa_{x^a}-g(x)\pa_s)\svop(\za'(x)s+\zb'(x))
\svop(\za(x)s+\zb(x))}$. It is obvious that the commutator bracket
of first-order differential operators induces a Lie algebroid
structure on $\sR\Z$ with the anchor
$(X\svop\zf\svop\zc)^0={\oX}$. In local affine coordinates:
\beas &[f_a(x)\pa_{x^a}-(\za(x)s+\zb(x))\pa_s+\zg(x),f'_a(x)\pa_{x^a}-
(\za'(x)s+\zb'(x))\pa_s+\zg'(x)]=\\
&(f_b(x)\frac{\pa f'_a}{\pa x^b}(x)-f'_b(x)\frac{\pa f_a}{\pa
x^b}(x))\pa_{x^a}-((f_a(x)\frac{\pa \za'}{\pa x^a}(x)-f'_a(x)\frac{\pa
\za}{\pa x^a}(x))s\\&+f_a(x)\frac{\pa \zb'}{\pa x^a}(x)-f'_a(x)\frac{\pa
\zb}{\pa x^a}(x)+\za(x)\zb'(x)-\za'(x)\zb(x))\pa_s\\
&+ (f_a(x)\frac{\pa \zg'}{\pa x^a}(x)-f'_a(x)\frac{\pa \zg}{\pa
x^a}(x))
\eeas
and
$$(f_a(x)\pa_{x^a}-(\za(x)s+\zb(x))\pa_s+\zg(x))^0=f_a(x)\pa_{x^a}.
$$
Writing $X=f_a(x)\pa_{x^a}$ and representing
$f_a(x)\pa_{x^a}-(\za(x)s+\zb(x))\pa_s+\zg(x)$ by
$(X,\za,\zb,\zg)$ we can write shortly
\be\label{222}
[(X,\za,\zb,\zg),(X',\za',\zb',\zg')]=([X,X'],X(\za')-X'(\za),X(\zb')-X'(\zb)+\za\zb'-\za'\zb,
X(\zg')-X'(\zg)).
\ee
Note that the distinguished sections $X_{{\sR\Z}}=-\pa_s$ and
$I_{{\sR\Z}}=1$ are in this Lie algebroid {\it ideal sections},
i.e. these sections are nowhere-vanishing and the sections of the
1-dimensional subbundles generated by $X_{{\sR\Z}}$ and
$I_{{\sR\Z}}$ are Lie ideals with respect to the Lie algebroid
bracket. A special vector bundle $(E,X_0)$ equipped with a Lie
algebroid structure such that $X_0$ is an ideal section we call an
{\it ideal-special Lie algebroid}. An ideal-special Lie algebroid
for which $X_0$ is a central section, i.e. $X_0$ commutes with any
section with respect to the Lie algebroid bracket, we call a {\it
special Lie algebroid}. It is easy to see that ideal-sections
define canonically 1-cocycles for the corresponding Lie
algebroids.

\begin{proposition} If $X_0$ is an ideal section of a Lie algebroid on the
vector bundle $E$ of rank $>1$ over $M$, then there is a closed
`1-form' $\zvf_{X_0}\in\Sec(E^\s)$ such that $[Y,X_0]=\lan
Y,\zvf_{X_0}\ran X_0$.
\end{proposition}
\begin{pf} Since $[X_0,fY]=f[X_0,Y]+\zr(X_0)(f)Y$ and $X_0$ generates a
Lie ideal, the anchor $\zr(X_0)$ vanishes if only $rank(E)>1$.
Thus $[Y,X_0]=\zF(Y)X_0$ for certain function $\zF(Y)$ which
linearly depends on $Y\in\Sec(E)$ and $\zF(fY)=f\zF(Y)$, so
$\zF(Y)=\lan Y,\zvf\ran$ for certain $\zvf\in\Sec(E^\s)$. The
`1-form' $\zvf$ is closed with respect to the Lie algebroid de
Rham differential, since, due to the Jacobi identity,
$\zF([Y_1,Y_2])=\zr(Y_1)(\zF(Y_2))-\zr(Y_2)(\zF(Y_1))$.
\end{pf}

\medskip\noindent
For the Lie algebroid ${\sR\Z}$ denote $\zvf_{X_{{\sR\Z}}}$ by $\zvf^0$.
In local affine coordinates, $\lan R,\zvf^0\ran=\za$ for
$R=f_a(x)\pa_{x^a}-(\za(x)s+\zb(x))\pa_s+\zg(x)$, so $\zvf^0$ is
nowhere-vanishing.  There is another canonical nowhere-vanishing closed
`1-form' $\zvf^1$ on ${\sR\Z}$ induced by the decomposition
$\de_{ah}(\Z)=\X_{ah}(\Z)\op Bas(\Z)$, namely $\lan\zvf^1,R\ran=\zg$. Note
that the form $\zvf_{I_{\sR\Z}}$ is identically zero.

The bundles $\wt{\sT}\Z$ and $\breve{\sL}\Z$ are subbundles of
$\sR\Z$ characterized by $\zvf^0=\zvf^1=0$ and $\zvf^1=0$,
respectively. On the level of realizations we have
$\wt{\sT}\Z=\wt{\sT}\Z\svop_M\bI \subset\breve{\sL}\Z$ and
$\breve{\sL}\Z=\breve{\sL}\Z\svop_M\bI\subset\sR\Z$. Of course,
$\wt{\sT}\Z$ and $\breve{\sL}\Z$ are Lie subalgebroids of $\sR\Z$
in every natural sense. Thus we have the chain
$\wt{\sT}\Z\subset\breve{\sL}\Z\subset\sR\Z$ of Lie algebroids
over $M$, canonically associated with $\Z$, whose Lie algebras of
sections are $\wt{\X}(\Z)$, $\X_{ah}(\Z)$, and $\de_{ah}(\Z)$,
respectively. The bundle $\Z^\dag$ is the kernel of the anchor map
in $\breve{\sL}\Z$, so $\Z^\dag$ is canonically a Lie algebroid
with the trivial anchor. It is easy to see that the Lie algebroid
bracket on $\Sec(\Z^\dag)=\Aff(\Z)$ is given by
\be\label{boz}
[\zf,\zc]=\{\zf,\zc\}_\Z.
\ee

\medskip\noindent
{\bf Remark.} The embedding of $\de_{ah}(\Z)$ into $\de(\Z)$
corresponds also to a Lie algebroid morphism from $\sR\Z$ into
$\sL\Z$. This morphism, however, is of a different kind than
morphism which are considered usually and which are associated
with the standard morphisms of vector bundles, and it is
represented by a relation, not a map. This kind of morphisms is
the Lie algebroid version of the Zakrzewski's morphisms of
groupoids (see \cite{Zak}). The Zakrzewski's morphisms of
groupoids lead to satisfactory functors into $C^\s$-algebras (cf.
\cite{Sta}).

\medskip
We can embed $\wt{\sT}\Z\svop_M\Z^\dag$ into
$\wt{\sT}\Z\svop_M\Z^\dag\svop_M\Z^\dag\simeq\sR\Z$ putting $\bI$
not on the third place but on the second. The resulting subbundle
of $\sR\Z$ we will denote $\wt{\sL}\Z$. It can be described as the
one determined by the equation $\zvf^1-\zvf^0=0$ and therefore it
is also a Lie subalgebroid of $\sR\Z$ like every kernel of a
closed nowhere-vanishing 1-form. The induced Lie algebroid
structure on $\wt{\sT}\Z\svop_M\Z^\dag$ reads
$$[X\svop\zf,X'\svop\zf']=[X,X']\svop(X(\zf')-X'(\zf)+\{\zf,\zf'\}_\Z)$$
and it is the same as the one obtained from the identification of
$\wt{\sT}\Z\svop_M\Z^\dag$ with $\breve{\sL}\Z$. In other words,
$\breve{\sL}\Z$ and $\wt{\sL}\Z$ are isomorphic Lie algebroids
differently placed in $\sR\Z$. The sections of $\wt{\sL}\Z$ are
first-order operators on $\Z$ having in local affine coordinates
the form
$$
{D}=f_a(x)\pa_{x^a}-(\za(x)s+\zb(x))\pa_s+\za(x).$$ The natural
isomorphism with $\breve{\sL}\Z$ is just the restriction of the
anchor map on $\sL\Z=\sT\Z\op\R$, i.e. $D\mapsto\oD$, where
$$
\oD=f_a(x)\pa_{x^a}-(\za(x)s+\zb(x))\pa_s.$$

\section{Affine derivations and affine first-order differential
operators}

Let us fix an AV-bundle $\zz:\Z\ra M$. In the standard
differential geometry the phase and the contact bundles $\sT^\s M$
and $\sT^\s M\op\R$ are representing objects for derivations and
linear first-order differential operators on $C^\infty(M)$, i.e.
on sections of the trivial vector bundle $M\ti\R$. By analogy, in
AV-differential geometry by the {\it bundle of affine derivations}
on $\Z$ (resp., the {\it bundle of affine first-order differential
operators} on $\Z$) with values in an affine bundle $A$ we
understand the affine bundle $\Aff_M(\sP\Z;A)$ (resp.,
$\Aff_M(\sC\Z;A)$). Thus the affine space $\ader(\Z;A)$ of {\it
affine derivations} (resp., the space $\ado(\Z;A)$ of {\it affine
first-order differential operators}) on $\Z$ with values in $A$ is
the space of sections of this bundle. We have an obvious action
$\zs\mapsto D(\zs)$ of $\sD\in\Aff(\sP\Z;A)$ (resp.
$\sD\in\Aff(\sC\Z;A)$) on sections $\zs$ of $\Z$ by
$D(\zs)=\sD(\bd\zs)$ (resp., $D(\zs)=\sD(\xc\zs)$). In the case
$A=M\ti\R$ we speak just about affine derivations (resp., affine
first-order operators) on $\Z$ and denote the (linear) space
$\ader(\Z;\R)=\Sec(\Aff_M(\sP\Z;\R))=\Sec(\sP\Z^\dag)$ (resp.
$\ado(\Z;\R)=\Sec(\Aff_M(\sC\Z;\R))=\Sec(\sC\Z^\dag)$) simply by
$\ader(\Z)$ (resp., $\ado(\Z)$). It is obvious by definition that
the linear parts of affine derivations (resp. differential
operators) are true derivations (resp. differential operators) on
$C^\infty(M)$. It is also clear that these concepts can be
extended naturally to a concept of a differential operator of
arbitrary order. In this sense, the affine space $\adoo(\Z;A)$ of
{\it affine differential operators of order 0} on $\Z$ with values
in $A$ is the space of sections of $\Aff_M(\Z,A)$, so the
differential operators of order 0 with values in $\R$ are sections
of $\Z^\dag$.

To understand better the structure of the bundles of derivations
and first-order differential operators let interpret them as
certain bundles constructed out of $\Z$ in the way in which
derivations of $C^\infty(M)$ are interpreted as vector fields,
i.e. sections of $\sT M$. Given an AV-bundle $\Z$ let us consider
the cotangent bundle $\sT^\s\Z$. The $(\R,+)$-action $\zvf$ on
$\Z$ can be lifted to an $(\R,+)$-action $\zvf^\s$ on $\sT^\s\Z$,
$(\zvf^\s)_r=(\zvf_{-r})^\s$. The fundamental vector field of this
action we denote by $X_{\sT^\s\Z}$. The orbits $[\za_{z_m}]$ of
this action form a vector bundle over $M$ which we denote by
$\widetilde{\sT}^\s\Z$. The sections of $\widetilde{\sT}^\s\Z$ are
represented by 1-forms on $\Z$, invariant with respect to
$X_{\sT^\s\Z}$. Moreover, there is a canonical decomposition
$\sT^\s\Z=\widetilde{\sT}^\s\Z\ti_M\Z$ given by
    \be\label{ts}\za_{z_m}\mapsto([\za_{z_m}],z_m)\ee
    which shows that
$\sT^\s\Z$ is canonically an affine bundle over $M$ with respect to the
projection $\zz\circ\zp_\Z$. This is a special affine bundle modelled on
$\widetilde{\sT}^\s\Z\ti{\bI}$. In local coordinates $(x^a,s)$ on $\Z$ and
the adapted coordinates $(x^a,s,p_a,p)$ on $\sT^\s\Z$, the lifted action
reads $(\zvf^\s)_r(x^a,s,p_a,p)=(x^a,s+r,p_a,p)$ and $X_{\sT^\s\Z}=-\pa_s$.
Hence, $(x^a,p_a,p)$ represent coordinates on $\widetilde{\sT}^\s\Z$ and the
section $p_a=p_a(x)$, $p=p(x)$, represents the invariant 1-form $p_a(x)\xd
x^a+p(x)\xd s$ on $\Z$. The affine phase bundle $\sP\Z$ can be identified
with the affine subbundle of $\widetilde{\sT}^\s\Z$ in obvious way:
$$\sP\Z=\{[\za_{z_m}]\in\widetilde{\sT}^\s\Z:
\lan\za_{z_m},X_\Z(z_m)\ran_{z_m}=1\}.
$$
Hence, $\widehat{\sP\Z}\simeq\widetilde{\sT}^\s\Z$. The contact
bundle $\sC\Z$ is an affine subbundle of
$\sT^\s\Z=\widetilde{\sT}^\s\Z\ati_M\Z$ being the affine product
$\sP\Z\ati_M\Z$.

We can do a similar procedure with the tangent bundle and obtain
$$\sT\Z=\widetilde{\sT}\Z\ti_M{\Z}.$$
The vector field $X_\Z$ is invariant, so it serves as a
distinguished section of $\widetilde{\sT}\Z$. Thus $\wt{\sT}\Z$ is
canonically a special vector bundle. Since $\widetilde{\sT}\Z$ is
dual to $\widetilde{\sT}^\s\Z$, it is obvious that
$\wt{\sT}\Z^\ddag=\sP\Z$ (or, equivalently that
$(\sP\Z)^\dag=\wt{\sT}\Z$), since sections of $\sP\Z$ are
considered as invariant 1-forms $\zn$ on $\Z$ such that
$\zn(X_\Z)=1$. Hence, $\ader(\Z)=\Sec(\wt{\sT}\Z)=\wt{\X}(\Z)$ is
the space of invariant vector fields $X$ on $\Z$ and their action
on sections $\zs$ of $\Z$ is given by $X(\zs)\circ\zz=X(\bF_\zs)$.
In local affine coordinates $(x^a,s)$ on $\Z$ for which
$X_\Z=-\pa_s$, we can write $X=f_a(x)\pa_{x^a}+g(x)X_\Z$, so that
$$
X(\zs)\circ\zz=(f_a(x)\pa_{x^a}-g(x)\pa_s)(\zs(x)-s)=
f_a(x)\frac{\pa\zs}{\pa x^a}(x)+g(x).
$$
We will use the natural convention and denote the pull-back
$f\circ\zz$ of a function $f\in C^\infty(M)$ to a basic function
on $\Z$ by $\bF_f$. With this convention we can simply write
$\bF_{X(\zs)}=X(\bF_\zs)$.

According to Theorem \ref{cc} (3), $\sC\Z^\dag$ equals
$$(\sP\Z)^\dag\svop_M\Z^\dag=\widetilde{\sT}\Z\svop_M\Z^\dag,$$ so
$\ado(\Z)=\Sec(\wt{\sT}\Z\svop_M\Z^\dag)$. The section
$D=X\svop\zf\in\Sec(\widetilde{\sT}\Z\svop_M\Z^\dag)$ acts on
$\zs\in\Sec(\Z)$ by $X(\zs)=X(\zs)+\zf\circ\zs$. We will identify
this bundle with $\wt{\sL}\Z$, since we can interpret this action
by ${D(\zs)}\circ\zz=(X+\bD_\zf)(\bF_\zs)$. In local affine
coordinates $D$ has the form
$$ D=f_a(x)\pa_{x^a}-(\za(x)s+\zb(x))\pa_s+\za(x)
$$
and its action on sections of $\Z$ reads
$$
\bF_{D(\zs)}=D(\zs)\circ\zz=D(\zs(x)-s)=f_a(x)\frac{\pa\zs}
{\pa{x^a}}(x)+\za(x)\zs(x)+\zb(x).
$$

In view of Corollary \ref{L} (Section 7),
$$\ader(\Z;\Z)=\Sec(\widetilde{\sT}\Z\bt_M\Z)=(\widetilde{\sT}\Z\ati_M\Z)/
\lan X_\Z-1_M\ran.
$$
Section $X\bt\zs'\in\Sec(\widetilde{\sT}\Z\bt_M\Z)$ acts on
$\zs\in\Sec(\Z)$ by $(X\bt\zs')(\zs)=X(\zs)+\zs'$. The embedding
$\bF$ of $\Z$ into $\Z^\dag$ induces the obvious embedding of
$\widetilde{\sT}\Z\bt_M\Z$ as an affine hyperbundle $\ol{\sT}\Z$
in $\widetilde{\sT}\Z\svop_M\Z^\dag$. If we identify the last
bundle with $\breve{\sL}\Z$, then $\ol{\sT}\Z$ can be interpreted
as an affine hyperbundle in $\breve{\sL}\Z$ and its  sections can
be  interpreted as first-order differential operators on $\Z$ of
the local form
$$\ol{X}=f_a(x)\pa_{x^a}+(s-\zb(x))\pa_s.
$$
Their action on sections of $\Z$ is given by
$$\bF_{\ol{X}(\zs)}(x,s)=\ol{X}(\zs(x)-s)=f_a(x)\frac{\pa\zs}
{\pa{x^a}}(x)+\zb(x)-s,
$$
so
$${\ol{X}(\zs)}(x)=f_a(x)\frac{\pa\zs}
{\pa{x^a}}(x)+\zb(x).
$$
Similarly as above, we get
$$\Aff(\sC\Z;\Z)=(\sC\Z)^\dag\bt_M\Z=\wt{\sL}\Z\bt_M\Z,$$ so that
$$\ado(\Z;\Z)=\Sec(\wt{\sL}\Z\bt_M\Z).$$
An element $\ol{D}=D\bt\zs'\in\Sec(\wt{\sL}\Z\bt_M\Z)$ acts on
$\zs\in\Sec(\Z)$ by $\ol{D}(\zs)=D(\zs)+\zs'$. Again, the
embedding $\bF:\Z\ra\Z^\dag$ induces the obvious embedding of
$\wt{\sL}\Z\bt_M\Z$ as an affine hyperbundle of the vector bundle
$\wt{\sL}\Z\svop_M\Z^\dag=\sR\Z$. This bundle we will denote
shortly $\ol{\sL}\Z$, so that, with respect to this
identification, $\ado(\Z;\Z)$ is the space of sections of
$\ol{\sL}\Z$, i.e. the space of first-order differential operators
on $\Z$ of the local form
$$\ol{D}=f_a(x)\pa_{x^a}-((\za(x)-1)s+\zb(x))\pa_s+\za(x).
$$
Then, $$\bF_{\ol{D}(\zs)}(x,s)=f_a(x)\frac{\pa\zs}
{\pa{x^a}}(x)+\za(x)\zs(x)+\zb(x)-s,
$$
so
$${\ol{D}(\zs)}(x)=f_a(x)\frac{\pa\zs}
{\pa{x^a}}(x)+\za(x)\zs(x)+\zb(x).
$$

We can summarize all these observations as follows.
\begin{theorem}\label{dob} Let $\Z$ be an AV-bundle over $M$.
There are subbundles of $\sR\Z$: vector subbundles  $\wt{\sT}\Z$,
$\wt{\sL}\Z$, $\breve{\sL}\Z$ and affine subbundles: $\ol{\sT}\Z$
modelled on $\wt{\sT}\Z$ and $\ol{\sL}\Z$ modelled on
$\wt{\sL}\Z$, characterized by $\zvf^0=\zvf^1=0$,
$\zvf^1-\zvf^0=0$, $\zvf^1=0$, $\zvf^0=1$ and $\zvf^1=0$,
$\zvf^1-\zvf^0=1$ , respectively, such that
\begin{description}
\item{(a)}
$\ader(\Z)=\Sec(\wt{\sT}\Z)=\wt{\X}(\Z)$,
\item{(b)}
$\ader(\Z;\Z)=\Sec(\ol{\sT}\Z)$,
\item{(c)}
$\ado(\Z)=\Sec(\wt{\sL}\Z)$,
\item{(d)}
$\ado(\Z;\Z)=\Sec(\ol{\sL}\Z)$,
\end{description}
and the action $\zs\mapsto D(\zs)$ of sections $D$ of these bundles,
regarded as elements of $\Sec(\sR\Z)=\de_{ah}(\Z)$, on sections $\zs$ of
$\Z$ is given by
$$\bF_{D(\zs)}=D(\bF_\zs).$$
\end{theorem}

\medskip\noindent
{\bf Remark.} Of course, the vector bundle $\Z^\dag$ (whose
sections represent $\adoo(\Z)$) and the affine bundle
$\Z^\dag\bt_M\Z$ (whose sections represent $\adoo(\Z;\Z)$) are
also subundles of $\sR\Z$ contained in the kernel of the anchor
map.

\section{Canonical Lie affgebroids associated
with AV-bundles}

In the standard differential geometry the canonical Lie algebroid
associated with a manifold $M$, or better, with the trivial bundle
$M\ti\R$, is $\sT M$. With an AV-bundle $\Z$ we have associated
the bundle $\wt{\sT}\Z$. Sections of $\wt{\sT}\Z$ are interpreted
as affine derivations on sections of $\Z$. The bundle $\wt{\sT}\Z$
carries a canonical  Lie algebroid structure like every Atiah
bundle of a principal bundle. The Lie bracket is inherited from
$\de(\Z)$. The bracket can be also described in terms of affine
derivations:
$$[X,X']=X_\sv\circ X'-(X')_\sv\circ X,$$
where $X_\sv$ is the vector part of the affine derivation
$X:\Sec(\Z)\ra C^\infty(M)$ (which represents also the anchor of
$X$).

Similarly, the bundle $\wt{\sL}\Z$ is also canonically a Lie algebroid
with similarly defined bracket
$$[D,D']=D_\sv\circ D'-(D')_\sv\circ D,$$
where $D_\sv:C^\infty(M)\ra C^\infty(M)$ is the vector part of
$D\in\ado(\Z)$.

Recall that the distinguished sections $X_{{\sR}\Z}=-\pa_s$ and
$I_{{\sR}\Z}=1$ are in the Lie algebroid  $\sR\Z$ {\it ideal
sections}, i.e. these sections are nowhere-vanishing and the
sections of the 1-dimensional subbundles generated by
$X_{{\sR}\Z}$ and $I_{{\sR}\Z}$ are Lie ideals with respect to the
Lie algebroid bracket. The closed `1-form' corresponding to
$X_{{\sR}\Z}$ we denote by $\zvf^0$.

The special affine bundles $\ol{\sT}\Z$ and $\ol{\sL}\Z$ also carry
canonical algebraic structures, represented by the commutators of their
sections regarded as affine maps $D:\Sec(\Z)\ra\Sec(\Z)$:
$$[D,D']=D\circ D'-D'\circ D.$$
These structures can be recognized as {\it Lie affgebroid} structures.
Recall (cf. \cite{GGU}) that an {\it affine Lie bracket} on  an  affine
space  $\A$ is a bi-affine map
   $$\left[\cdot,\cdot\right]\colon  \A\times\A\ra \sV(\A)$$
which is skew-symmetric:
$\left[\zs_1,\zs_2\right]=-\left[\zs_2,\zs_1\right]$ and satisfies
the Jacobi identity:
        $$\left[\zs_1,\left[\zs_2, \zs_3\right]\right]_\sv^2+\left[\zs_2,
\left[\zs_3, \zs_1\right] \right]_\sv^2+ \left[\zs_3,\left[\zs_1,
\zs_2\right]\right]_\sv^2=0,$$ where $[\cdot,\cdot]_\sv^2$ is the
affine-linear part of the biaffine bracket.
   An affine space equipped with an affine Lie bracket we shall call a
{\it Lie affgebra}. Note that the term {\it affine Lie algebra}
has been  already  used  for certain types of Kac-Moody algebras.

   If $A$ is an affine bundle over $M$ modelled on $\sV(A)$ then a {\it Lie
affgebroid structure} on $A$ is an affine Lie bracket on sections
of $A$ and a morphism $\zg\colon A\ra\sT M$ of affine bundles
(over the identity on $M$) such that $[\zs,\cdot]_\sv^2$ is a
quasi-derivation with the anchor $\zg(\zs)$, i.e.
      $$ \left[\zs,fX\right]_\sv^2=f\left[\zs,X\right]_\sv^2+\zg(\zs)(f)X$$
   for all $\zs\in\Sec(A)$, $X\in\Sec(\sV(A))$, $f\in C^\infty(M)$.

\bigskip\noindent{\bf Remark.}
The above definition is a slight generalization of the one
proposed in \cite{MMS1,MMS2} where the additional assumptions that
the base manifold $M$ is fibered over $\R$  and   that $\zg(\zs)$
are vector fields projectable onto $\frac{\partial}{\partial t}$
have been put. On the other hand, one can try to define Lie
affgebra as a skew-symmetric (in the affine sense) bracket
$[\cdot,\cdot]_a$ on $\Sec(A)$ with values in $\Sec(A)$ satisfying
the Jacobi identity of the form
$$\sum_{\zw\in
S_3}sgn(\zw)[\zs_{\zw(1)},[\zs_{\zw(2)},\zs_{\zw(3)}]_a]_a=0.$$ The l.h.s.
of the above equation is a vector combination of elements of $\Sec(A)$, so
the identity makes sense. The problem with such definition is that, as we
already know, any skew-symmetric operation on $\Sec(A)$ defines
automatically an element $\zs_0\in\Sec(A)$, $\zs_0=[\zs,\zs]$, and we get
such a bracket in the form $[\zs_1,\zs_2]_a=[\zs_1,\zs_2]+\zs_0$, where
$[\cdot,\cdot]$ is the Lie affgebra bracket in the version we started
with. Fixing $\zs_0$ is usually too much (we just get a trivialization of
the affine space) for applications and canonical examples, so we remain
with the weaker definition.

\medskip\noindent
{\bf Example 1.} Every AV-bundle $\Z$ carries a canonical Lie
affgebroid structure induced by the affine structure. The bracket
of sections $\zs,\zs'$ of $\Z$ is just $[\zs,\zs']=(\zs-\zs')$.

\bigskip\noindent
The following  fact  has  been  proved  in \cite{GGU}, Theorem 11.
\begin{theorem} A map $\left[\cdot,\cdot\right]\colon \Sec(A)\times
\Sec(A)\ra\Sec(\sV(A))$ is a Lie affgebroid bracket on an affine bundle
$A$ if and only if there is an extension of this map to a Lie algebroid
bracket $[\cdot,\cdot]^\we$ on $\wh{A}$ such that
\be\label{aff}[\Sec(\wh{A}),\Sec(\wh{A})]^\we\subset\Sec(\sV(A)).
\ee
Moreover, (\ref{aff}) is equivalent to the fact that
$1_A\in\Sec(A^\dag)=\Sec(\wh{A}^\s)$ is a closed one-form.
\end{theorem}
The Lie algebroid $(\wh{A},[\cdot,\cdot]^\we)$ is uniquely
determined by the Lie affgebroid $({A},[\cdot,\cdot])$ and we will
call it the {\it Lie algebroid hull} of $({A},[\cdot,\cdot])$.

\medskip\noindent
{\bf Example 2.} The Lie affgebroid bracket on $\Z$ from the previous
example extends to a Lie algebroid bracket on $\wh{\Z}$. This bracket can
be expressed by means of the canonical Jacobi bracket on $\Z$ by
$\bF_{[u,u']}=\{\bF_u,\bF_{u'}\}_\Z$.

\bigskip\noindent

Since the affine subbundles $\ol{\sT}\Z$ and $\ol{\sL}\Z$ in the
Lie algebroids $\breve{\sL}\Z$ and $\wt{\sR}\Z$ are defined as the
1-level sets of $\zvf^0$ and $\zvf^1-\zvf^0$, respectively, we get
the following (cf. \cite{GGU}).

\begin{theorem} The special affine bundles $\ol{\sT}\Z$ and
$\ol{\sL}\Z$ carry canonical Lie affgebroid structures for which
the brackets are the commutators in $\ader(\Z;\Z)$ and
$\ado(\Z;\Z)$, respectively. The Lie affgebroid hulls of
$\ol{\sT}\Z$ and $\ol{\sL}\Z$ are $\breve{\sL}\Z$ and
$\wt{\sR}\Z$, respectively.
\end{theorem}

\section{Aff-Poisson and aff-Jacobi brackets}
The idea of an affine analog of a Poisson bracket goes back to
\cite{Ur1} but we will mainly follow the picture described in
\cite{GGU}.

Let $\Z$ be an AV-bundle over $M$. An affine Lie bracket on
$\sss{Sec}(\Z)$
        $$      \{\cdot,\cdot\} \colon \sss{Sec}(\Z)\times \sss{Sec}(\Z) \rightarrow
C^\infty(M)$$
        is called an {\it aff-Poisson} (resp.{\it aff-Jacobi}) {\it
bracket} if
        $$ \{\zs,\cdot \} \colon \sss{Sec}(\Z) \rightarrow C^\infty(M)$$
        is an affine derivation (resp. an affine first-order differential
operator) for every $\zs\in \sss{Sec}(\Z)$.

\medskip\noindent We use  the  term {\it aff-Poisson}, since  {\it  affine
Poisson structure}  has already  a different meaning in the literature.

\medskip\noindent
{\bf Example 3.} Every AV-bundle $\Z$ carries a canonical
aff-Jacobi bracket determined by the affine structure:
\be\label{caj}\{\zs,\zs'\}=\zs-\zs'.
\ee

\begin{theorem}\cite{GGU}\label{apj}
        For every aff-Poisson (resp. aff-Jacobi) bracket
        $$\{\cdot,\cdot\} \colon \sss{Sec}(\Z)\times \sss{Sec}(\Z) \rightarrow
C^\infty(M)$$
        its vector part
        $$\{\cdot,\cdot\}_\sv \colon C^\infty(M)\times C^\infty(M) \rightarrow
C^\infty(M)$$
        is a Poisson (resp. Jacobi) bracket.  Moreover,
         $$ \{\zs,\cdot \}_\sv^2 \colon C^\infty(M) \rightarrow C^\infty(M)
                                         $$
        is a derivation (resp. first-order differential operator) for every
section $\zs\in \sss{Sec}(\Z)$, which is simultaneously a derivation of
the bracket  $\{\cdot,\cdot\}_\sv$.  Conversely,  if we  have a Poisson
(resp. Jacobi) bracket $\{\cdot,\cdot\}_0$  on $C^\infty(M)$ and  a
derivation  (resp.  a first-order differential operator)
          $$ D \colon C^\infty(M) \rightarrow C^\infty(M)
                                         $$
which is simultaneously a derivation of the bracket $\{\cdot,\cdot\}_0$,
then there is a unique aff-Poisson (resp. aff-Jacobi)  bracket
$\{\cdot,\cdot\}$ on  $\sss{Sec}(\Z)$ such that
$\{\cdot,\cdot\}_0=\{\cdot,\cdot\}_\sv$ and $D=\{\zs,\cdot\}_\sv^2$ for a
chosen section $\zs\in \sss{Sec}(\Z)$.

Using a section $\zs_0$ to identify $\sss{Sec}(\Z)$ with $C^\infty(M)$, we
get that the aff-Poisson (resp. aff-Jacobi) bracket on $\sss{Sec}(\Z)$ has
the form
                $$ \{\zs,\zs'\} = D(\zs'-\zs) + \{\zs,\zs'\}_\sv,$$
        where $D$ is a vector field (resp first-order differential operator)
which  is  a  derivation  of  the   Poisson   (resp. Jacobi) bracket
$\{\cdot,\cdot\}_\sv$.
          \end{theorem}

\medskip\noindent
{\bf Example 4.} Every Poisson (resp., Jacobi) bracket
$\{\cdot,\cdot\}_M$ on $C^\infty(M)$ can be interpreted as an
aff-Poisson (resp., aff-Jacobi) bracket $\{\cdot,\cdot\}$ on
sections of the trivial AV-affine bundle $M\ti\bI$. In this case
the trivialization is canonical, $D=0$ and
$\{\cdot,\cdot\}=\{\cdot,\cdot\}_M$.

\begin{theorem}\label{m1} Let $\zz:\Z\ra M$ be an AV-bundle. Then,
\begin{description}
\item{(1)} there is a one-to-one correspondence between aff-Poisson
brackets $\{\cdot,\cdot\}_{aP}$ on $\Sec(\Z)$ and Poisson brackets
$\{\cdot,\cdot\}_{\zP}$ on $C^\infty(\Z)$ which are $X_\Z$-invariant, i.e.
which are associated with Poisson tensors $\zP$ on $\Z$ such that
$\Ll_{X_\Z}\zP=0$. This correspondence is determined by
\be\label{ap1}\{\zs,\zs'\}_{aP}\circ\zz=\{\bF_{\zs},\bF_{\zs'}\}_\zP;
\ee
\item{(2)} there is a one-to-one correspondence between aff-Jacobi
brackets $\{\cdot,\cdot\}_{aJ}$ on $\Sec(\Z)$ and Jacobi brackets
$\{\cdot,\cdot\}_{J}$ on $C^\infty(\Z)$ which are associated with Jacobi
structures $J=(\zP,\zG)$ on $\Z$ such that $\Ll_{X_\Z}\zG=0$ and
$\Ll_{X_\Z}\zP=\zG\we X_\Z$. This correspondence is determined by
\be\label{apj1} \{\zs,\zs'\}_{aJ}\circ\zz=\{\bF_{\zs},\bF_{\zs'}\}_J; \ee
\end{description}
\end{theorem}
\begin{pf} We will prove only part (2). The proof of (1) is analogous
but easier. Since all objects are local over $M$, we can use local
affine coordinates $(x^a,s)$ on $\Z$ in which $X_\Z=-\pa_s$ and
identify sections $\zs$ of $\Z$ with functions $\zs(x)$, so that
$\bF_\zs(x,s)=\zs(x)-s$. We will identify functions on $M$ with
basic functions on $\Z$. Assume first that $\{\cdot,\cdot\}_{aJ}$
is an aff-Jacobi bracket on $\Sec(\Z)$. According to Theorem
\ref{apj} there is a Jacobi structure $J_0=(\zP_0,\zG_0)$ on $M$
and a first-order differential operator $D={\oD}+f$ on $M$ such
that $\{\zs,\zs'\}_{aJ}=D(\zs-\zs')+\{\zs,\zs'\}_{J_0}$. The
equation (\ref{apj1}) can be rewritten in the form
\beas&\zP_0(\zs,\zs')+\zs\zG_0(\zs')-\zs'\zG_0(\zs)+{\oD}(\zs-\zs')+f(\zs-\zs')=\\
&\zP(\zs-s,\zs'-s)+(\zs-s)\zG(\zs'-s)-(\zs'-s)\zG(\zs-s), \eeas which has
a unique solution $J=(\zP,\zG)$, namely
$$
\zP=\zP_0+\pa_s\we({\oD}-s\zG_0),\quad \zG=\zG_0-f\pa_s.
$$
It is easy to
see that the Jacobi identity for $\{\cdot,\cdot\}_{aJ}$ implies the Jacobi
identity for $\{\cdot,\cdot\}_{J}$ and that this solution has the required
properties with respect to $\Ll_{X_\Z}$.

Conversely, assume that $J$ is a Jacobi structure on $\Z$. The
conditions $\Ll_{X_\Z}\zG=0$ and $\Ll_{X_\Z}\zP=\zG\we X_\Z$ imply
that $\Ll_{X_\Z}(\{(\zs-s),(\zs'-s)\}_J)=0$, i.e.
$\{(\zs-s),(\zs'-s)\}_J$ is a basic function, so that (\ref{apj})
defines a bracket on $\Sec(\Z)$. It is easy to see that this
bracket is an aff-Jacobi bracket.
\end{pf}

\medskip
Note, that every skew-symmetric affine bracket $\{\cdot,\cdot\}$
is uniquely determined by  $\{\zs,\cdot \}_\sv^2$, namely
    \be \{\zs,\zs'\} =  \{\zs,\zs' -\zs \}_\sv^2 .\ee
For an aff-Poisson bracket on sections of $\Z$ the  mapping
$f\mapsto \{\zs,f \}_\sv^2$ is a derivation of the algebra
$C^\infty(M)$, hence a vector field on $M$. We denote it by
$X_\zs$ and we call it the {\it hamiltonian vector field} of
$\zs$.

\medskip\noindent
{\bf Example 5.} The canonical Poisson structure $\zP$ on the
cotangent bundle $\sT^\s M$ is invariant with respect to
translation by the vertical lift $\za_{\sT^\s M}$ of any closed
1-form $\za\in\Sec(\sT^\s M)$. If $\za$ is nowhere-vanishing, we
can consider the corresponding AV-bundle $\Z=\bAP(\sT^\s M)$ for
which $X_\Z$ is the vertical lift of $\za$, i.e. $X_\Z=\za_{\sT^\s
M}$. Hence, the AV-bundle $\Z$, i.e. $\zz:\sT^\s M\ra\sT^\s
M/\la\za\ran$, carries a canonical aff-Poisson structure with the
bracket (\ref{ap1}). Since, for any section $\zs$ of $\Z$ and for
any function $f$ on $\sT^\s M/\la\za\ran$, we have
$(\{\zs,f\}_{aP})^2_\sv\circ\zz=\{\bF_\zs,f\circ\zz\}_\zP$, the
hamiltonian vector field $X_\zs$ on $\sT^\s M/\la\za\ran$ induced
by the section $\zs$ is the projection of the hamiltonian vector
field on $\sT^\s M$ induced by $\bF_\zs$.

\medskip
In the theory of Lie algebroids it is well-known that a Lie
algebroid brackets $[\cdot,\cdot]$ on the vector bundle $E$ are in
a one-to-one correspondence with linear Poisson brackets
$\{\cdot,\cdot\}$ on $E^\s$. Linearity of the bracket means that
the bracket of linear functions is a linear function and the
correspondence is described by
$$\{\zi_{E^\s}(X_1),\zi_{E^\s}(X_2)\}=\zi_{E^\s}([X_1,X_2]), $$ where
$\zi_{E^\s}(X)$ denotes the linear function on $E^\s$ associated
canonically with $X\in\Sec(E)$. In \cite{GIM} it has been shown that this
correspondence can be extended to a one-to-one correspondence between Lie
algebroid brackets on $E$ and affine Jacobi brackets (bracket of affine
functions is an affine function) on an arbitrary affine hyperbundle $A$ of
$E^\s$. In the case of special affine bundles we have an analogous
correspondence which refers to Theorem \ref{oto}.

Let $\bA=(A,v^0)$ be a special affine bundle over $M$. There is an
obvious identification of a section $X$ of $\sV(\bA)$ with a
linear function $\zi_{\bA^\dag}(X)$  on $\bA^\dag$ and an affine
function $\zi_{\bA^\#}(X)$ on $\bA^\#$ which are invariant with
respect to translation by $1_\bA$, so they are pull-backs of a
certain linear function $\zi_X^\dag$ and an affine function
$\zi_X^\#$ on $\bA^\dag/\langle 1_\bA\rangle$ and $\bA^\#/\langle
1_\bA\rangle$, respectively.

\begin{theorem}\label{c3}\cite{GGU} There is a one-to-one correspondence between
Lie affgebroid brackets $[\cdot,\cdot]_\bA$ on a special affine bundle
$\bA$ and
\begin{description}
\item{(1)} linear aff-Poisson brackets $\{\cdot,\cdot\}_{\bA^\dag}$ on the
AV-bundle ${\bAP(\bA^\dag)}$, i.e. on $ \zr^\dag\colon \bA^\dag
\rightarrow \bA^\dag/\langle 1_\bA\rangle$, determined by
\be\label{c1}\zi^\dag_{[a,a']_\bA}=
\{\wh{\zs}_{a},\wh{\zs}_{a'}\}_{\bA^\dag};\ee
\item{(2)} affine aff-Jacobi brackets $\{\cdot,\cdot\}_{\bA^\#}$ on the
AV-bundle ${\bAP(\bA^\#)}$, i.e. on $ \zr\colon \bA^\# \rightarrow
\bA^\#/\langle 1_\bA\rangle$, determined by
\be\label{c2}\zi^\#_{[a,a']_\bA}=\{{\zs}_{a},{\zs}_{a'}\}_{\bA^\#}.\ee
This aff-Jacobi bracket is aff-Poisson if and only if the section
$v_0$ is central in the Lie algebroid hull $\wh{\bA}$ of $\bA$.
\end{description}
\end{theorem}

\medskip\noindent
{\bf Remark.} Here we call an aff-Poisson (resp., aff-Jacobi)
structure linear (resp., affine) if the bracket of linear sections
of $ \zr^\dag\colon \bA^\dag \rightarrow \bA^\dag/\langle
1_\bA\rangle$ (resp., the bracket of affine sections of $
\zr\colon \bA^\# \rightarrow \bA^\#/\langle 1_\bA\rangle$) is a
linear function on $\bA^\dag /\langle 1_\bA\rangle $ (resp., it is
an affine function on $\bA^\#/\langle 1_\bA\rangle$).

\medskip\noindent
\begin{pf} Part (1) has been proved in \cite{GGU}, Theorem 19, so
$[\cdot,\cdot]_\bA$ induces certain aff-Poisson bracket
$\{\cdot,\cdot\}_{\bA^\dag}$. According to Theorem \ref{m1}, there
is a Poisson tensor $\zP$ on $\bA^\dag$ which corresponds to the
aff-Poisson bracket $\{\cdot,\cdot\}_{\bA^\dag}$. This tensor is,
clearly, linear and it is invariant with respect to the vector
field $X^\dag=(1_\bA)_{\bA^\dag}$ - the vertical lift of the
section $1_\bA$. Now, we use the result \cite{GIM}, Corollary 3.6,
which implies that there is a one-to-one correspondence between
linear Poisson brackets $\{\cdot,\cdot\}_\zP$ on $\bA^\dag$ and
affine Jacobi brackets $\{\cdot,\cdot\}_J$ on $\bA^\#$ such that
$$
\{u,v\}_\zP\mid_{\bA^\#}=\{ u\mid_{\bA^\#},v\mid_{\bA^\#}\}_J
$$
for all linear functions $u,v$ on $A^\dag$. The Jacobi structure
$J=(\zP_0,\zG_0)$ is the restriction to $\bA^\#$ of the Jacobi
structure $(\zP+\zG\we\zD_{\bA^\dag},\zG)$, where $\zG=\{
v^0,\cdot\}_\zP$ is the Hamiltonian vector field of the linear
function $\zi_{\bA^\dag}(v^0)\in\Sec(\sV(\bA))$ which defines the
affine hyperbundle $\bA^\#$ in $\bA^\dag$ and $\zD_{\bA^\dag}$ is
the Liouville (Euler) vector field on the vector bundle
$\bA^\dag$. The crucial point is that $X^\dag$ preserves $\zP$ if
and only if it preserves $\zG$ and
\be\label{78}
\Ll_{X^\dag}(\zP+\zG\we\zD_{\bA^\dag})=\zG\we X^\dag.
\ee
Indeed, since the vector field $X^\dag$ preserves $\zP$, and the function
$\zi_{\bA^\dag}(v^0)$ due to the fact that $X^\dag(\zi_{\bA^\dag}(v^0))$
is the pull-back of $\lan 1_\bA,v^0\ran=0$, it preserves also $\zG$, i.e.
$\Ll_{X^\dag}\zG=0$. Moreover, since $X^\dag$ is a vertical lift,
$\Ll_{X^\dag}\zD_{\bA^\dag}=X^\dag$. Thus we get (\ref{78}) and, due to
Theorem \ref{c3}(b), we get an aff-Jacobi bracket
$\{\cdot,\cdot\}_{\bA^\#}$ on ${\bAP(\bA^\#)}$. It is easy to see that it
satisfies (\ref{c2}). The converse is proved by a similar reasoning in
reversed order. Passing to the restrictions to $\bA^\#$ we get, in view of
(\ref{m1}), that $J$ corresponds to an aff-Jacobi bracket on sections of
$\zr$. Since $\bF^\dag_{\wh{\zs}_a}\mid_{\bA^\#}=\bF^\#_{\zs_a}$, the
theorem is proved.
\end{pf}

\medskip\noindent
{\bf Example 6.} The canonical Lie affgebroid structure on $\Z$
given by $[\zs,\zs']=(\zs'-\zs)$ induces an aff-Poisson structure
on the AV-bundle $\bAP(\Z^\dag)$ and an aff-Jacobi structure on
$\bAP(\Z^\#)=\Z$. The corresponding linear Poisson structure on
$\Z^\dag$ (resp., affine Jacobi structure on $\Z$) is
$\zP=\zD_{\Z^\dag}\we X_\Z$ (resp., $J=(0,X_\Z)$).

\medskip\noindent
{\bf Example 7.} Consider the Lie algebroid structure on
$\wt{\sT}\Z$ as a Lie affgebroid structure on the special affine
bundle (in fact, special vector bundle). The special affine dual
$(\wt{\sT}\Z)^\#$ is $\sP\Z\ati\bI$, so $\bAP((\wt{\sT}\Z)^\#)$ is
the trivial AV-bundle over $\sP\Z$. The corresponding aff-Jacobi
bracket on $\sP\Z\ati\bI$ is the aff-Poisson bracket induced from
the canonical Poisson structure on $\sP\Z$ associated with the
canonical symplectic structure.

\medskip\noindent
{\bf Example 8.} Recall that we have the identification
$\ol{\sT}\Z=\wt{\sT}\Z\bt_M\Z=(\sP\Z)^\dag\bt_M\Z$, so that, according to
Theorem \ref{bt}, $(\ol{\sT}\Z)^\#=\sP\Z\ati_M\Z^\#$. The corresponding
AV-bundle is $\zr:\sP\Z\ati_M\oZ^\#\ra\sP\Z$. But
$\sP\Z\ati_M\oZ^\#=\sP\Z\ati_M\Z=\sC\Z$, so that
${\bAP((\ol{\sT}\Z)^\#)}=\bAP(\sC\Z)$. The special affine bundle
$\ol{\sT}\Z$ is canonically a Lie affgebroid, so, due to the above
theorem, $\bAP(\sC\Z)$ is equipped with a canonical aff-Jacobi structure.
It is easy to guess that this is the structure corresponding, {\it via}
Theorem \ref{m1}, to the canonical Jacobi bracket (\ref{cjb}) on $\sC\Z$
associated with the canonical contact form which, in turn, is represented
by the affine Liouville 1-form. Indeed, let $(x^a,p_b,s)$ be standard
affine coordinates on $\sC\Z$ induced from the Darboux coordinates in
$\sT^\s M\ti\R$ and let $(x^a,f_b,\zb)$ be the coordinates in $\ol{\sT}\Z$
representing the vector field
$${X}=f_a(x)\pa_{x^a}+(s-\zb(x))\pa_s\in\X_{ah}(\Z).$$ The duality
between $\ol{\sT}\Z=(\sP\Z)^\dag\bt_M\Z$ and $\sC\Z=\sP\Z\ati_M\Z$ is
given by $\lan(x^a,f_b,\zb),(x^a,p_b,s)\ran_{as}=f_ap_a+\zb-s$ so that
$\zi^\#_{\sC\Z}(X)(x^a,p_b,s)=f_a(x)p_a+\zb(x)-s$ and
$\zs_X(x^a,p_b)=(x^a,p_b,f_a(x)p_a+\zb(x))$. Since the Lie affgebroid
bracket in $\ol{\sT}\Z$ reads
\beas[X,X']_{\ol{\sT}\Z}&=&[f_a(x)\pa_{x^a}+(s-\zb(x))\pa_s,
f'_b(x)\pa_{x^b}+(s-\zb'(x))\pa_s]_{\de(\Z)}\\
&=&(f_b(x)\frac{\pa f'_a}{\pa x^b}(x)-f'_b(x)\frac{\pa f_a}{\pa
x^b}(x))\pa_{x^a}-(f_a(x)\frac{\pa \zb'}{\pa x^a}(x)-f'_a(x)\frac{\pa
\zb}{\pa x^a}(x)+\zb(x)-\zb'(x))\pa_s,
\eeas
the corresponding Jacobi bracket on $\sC\Z$ is uniquely characterized by
\beas\{ f_a(x)p_a+\zb(x)-s,f'_b(x)p_b+\zb'(x)-s\}_J&=&
p_a(f_b(x)\frac{\pa f'_a}{\pa x^b}(x)-f'_b(x)\frac{\pa f_a}{\pa
x^b}(x))\\&&+(f_a(x)\frac{\pa \zb'}{\pa x^a}(x)-f'_a(x)\frac{\pa \zb}{\pa
x^a}(x)+\zb(x)-\zb'(x)).
\eeas
It is easy to check that this is exactly the Jacobi structure
$$
J=(\pa_{p_a}\we\pa_{x^a}+p_a\pa_{p_a}\we\pa_s,-\pa_s),$$ i.e. the Jacobi
structure of the contact 1-form $p_a\D x^a-\D s$.

\section{Aff-Poisson and aff-Jacobi (co)homology}

Let $\Z$ be an AV-bundle over $M$. It is obvious that affine
biderivations on $\Z$ are affine derivations on $\Z$ with values
in $\ader(\Z)$, i.e. sections of the bundle
$$\Aff_M(\sP\Z;\wt{\sT}\Z)=\Hom_M(\wh{\sP\Z};\wt{\sT}\Z))=
(\sP\Z)^\dag\ot_M\wt{\sT}\Z=\wt{\sT}\Z\ot_M\wt{\sT}\Z.
$$
In this picture, skew-symmetric affine biderivations are sections
of $\bigwedge^2\wt{\sT}\Z$. Similarly, affine first-order
bidifferential operators on $\Z$ are sections of
$\bigwedge^2\wt{\sL}\Z$. Since both, $\wt{\sT}\Z$ and
$\wt{\sL}\Z$, are Lie algebroids, there are the corresponding Lie
algebroid Schouten brackets $\lna\cdot,\cdot\rna_{\wt{\sT}\Z}$ and
$\lna\cdot,\cdot\rna_{\wt{\sL}\Z}$ on the Grassmann algebras
$\sA(\wt{\sT}\Z)=\bigoplus_n\sA^n(\wt{\sT}\Z)=\bigoplus_n\Sec(\bigwedge^n\wt{\sT}\Z)$
and $\sA(\wt{\sL}\Z)=\bigoplus_n\Sec(\bigwedge^n\wt{\sL}\Z)$ of
the vector bundles $\wt{\sT}\Z$ and $\wt{\sL}\Z$, respectively.
What will be crucial here is that the Lie algebroid $\wt{\sL}\Z$
possesses additionally a canonical closed `1-form' $\zvf^0$
inherited from $\sR\Z$ (in fact, $\zvf^0=\zvf^1$ on $\wt{\sL}\Z$)
which makes it into a {\it Jacobi algebroid} with the
Schouten-Jacobi bracket
$\lna\cdot,\cdot\rna_{\wt{\sL}\Z}^{\zvf^0}$.

\medskip\noindent
{\bf Remark.} The Jacobi algebroids have been introduced by
Iglesias and Marrero \cite{IM} under the name {\it generalized Lie
algebroids} and recognized as graded Jacobi brackets in
\cite{GM1,GM2}. For the definitions and details we refer to these
papers or to the article \cite{GU3} which contains a brief
introduction to the theory of Jacobi algebroids, the corresponding
lifts of tensors and {\it canonical structures}.

\begin{theorem}\
\begin{description}
\item{(a)} $\zL\in\Sec(\bigwedge^2\wt{\sT}\Z)$ represents an aff-Poisson
structure on $\Z$ if and only if $\lna\zL,\zL\rna_{\wt{\sT}\Z}=0$.
\item{(b)} $\cJ\in\Sec(\bigwedge^2\wt{\sL}\Z)$ represents an aff-Jacobi
structure on $\Z$ if and only if
$\lna\cJ,\cJ\rna_{\wt{\sL}\Z}^{\zvf^0}=0$.
\end{description}
In other words, aff-Poisson and aff-Jacobi structures are
canonical structures for the Lie algebroid $\wt{\sT}\Z$ and the
Jacobi algebroid $(\wt{\sL}\Z,\zvf^0)$, respectively.
\end{theorem}
\begin{pf} (a) We will use a trivialization of $\Z$ to identify
$\wt{\sT}\Z$ with $\sL M=\sT M\op\R$ and we will use the
expression $D=(X,\zb)\in\Sec(\sL M)=\X(M)\ti C^\infty(M)$ for
sections $D$ of $\wt{\sT}\Z$ (see the convention preceding
(\ref{222})). The action on $\zs\in\Sec(\Z)$, identified with
functions on $M$, reads $D(\zs)=X(\zs)+\zb$. With respect to this
identification, sections of $\wt{\sT}\Z$ commute exactly as
sections of $\sL M$, i.e. (cf. (\ref{222}))
$$[(X,\zb),(X',\zb')]_{\wt{\sT}\Z}=([X,X']_{\sT M},X(\zb')-X'(\zb)).$$
Since $\bigwedge^2\wt{\sT}\Z$ is identified with
$\bigwedge^2({\sT}M\op\R)$, elements
$\zL\in\Sec(\bigwedge^2\wt{\sT}\Z)$ are of the form
$\zL=\zL_0+X_\Z\we X_0$, where $X_\Z=(0,1)$,
$\zL_0\in\Sec(\bigwedge^2{\sT}M)$ is a bivector field on $M$ and
$X_0$ is a vector field on $M$. The bi-section $\zL$ induces the
bracket
$$\{\zs,\zs'\}_\zL=\{\zs,\zs'\}_{\zL_0}+X_0(\zs'-\zs).$$
In view of Theorem \ref{apj} this is an aff-Poisson bracket if and
only if $\zL_0$ is a Poisson tensor and $\lna
X_0,\zL_0\rna^{SN}=0$, where $\lna\cdot,\cdot\rna^{SN}$ is the
Schouten-Nijenhuis bracket, i.e. the Lie algebroid Schouten
bracket for $\sT M$. But these conditions are equivalent to
$\lna\zL,\zL\rna_{\wt{\sT}\Z}=0$. Indeed, since $X_\Z$ is a
central section,
\beas
\lna\zL,\zL\rna_{\wt{\sT}\Z}&=&\lna\zL_0,\zL_0\rna_{\wt{\sT}\Z}+
2\lna\zL_0,X_\Z\we X_0\rna_{\wt{\sT}\Z}+\lna X_\Z\we X_0 ,X_\Z\we
X_0 \rna_{\wt{\sT}\Z}\\
&=&\lna\zL_0,\zL_0\rna^{SN}-2X_\Z\we\lna\zL_0,X_0\rna^{SN},
\eeas
that vanishes exactly when $\lna\zL_0,\zL_0\rna^{SN}=0$ and
$\lna\zL_0,X_0\rna^{SN}=0$.

\medskip\noindent
(b) Similarly as above we use an identification $\wt{\sL}\Z\simeq
\sL\op\R$ and the expression $D=(X,\zb)\in\Sec(\sL
M\op_M\R)=\de(M)\ti C^\infty(M)$ for sections $D$ of $\wt{\sL}\Z$.
The action on $\zs\in\Sec(\Z)$ reads $D(\zs)=X(\zs)+\zb$. With
respect to this identification, sections of $\wt{\sL}\Z$ commute
like
$$[(X,\zb),(X',\zb')]_{\wt{\sL}\Z}=([X,X']_{\sL M},X(\zb')-X'(\zb)).$$
Elements $\cJ\in\Sec(\bigwedge^2\wt{\sL}\Z)$ are of the form
$\cJ=\cJ_0+X_\Z\we D_0$, where $X_\Z=(0,1)$,
$\cJ_0\in\Sec(\bigwedge^2{\sL}M)$ is a first-order bidifferential
operator on $M$ and $D_0$ is a first-order differential operator
on $M$. The Schouten-Jacobi bracket
$\lna\cdot,\cdot\rna_{\wt{\sL}\Z}^{\zvf^0}$ restricted to $\sL M$
gives the canonical Schouten-Jacobi bracket
$\lna\cdot,\cdot\rna_{{\sL}M}^{\zvf^0}$ on $\sA(\sL M)$ for which
canonical structures are Jacobi structures on $M$ (cf.
\cite{GM1}). The section $\cJ$ induces the bracket
$$\{\zs,\zs'\}_\cJ=\{\zs,\zs'\}_{\cJ_0}+D_0(\zs'-\zs).$$
This is an aff-Jacobi bracket if and only if $\cJ_0$ is a Jacobi
structure and $\lna D_0,\cJ_0\rna_{{\sL}M}^{\zvf^0}=0$ (Theorem
\ref{apj}). But these conditions are equivalent to
$\lna\cJ,\cJ\rna_{\wt{\sL}\Z}^{\zvf^0}=0$. Indeed, since $X_\Z$ is
an ideal section such that $[D,X_\Z]_{\wt{\sL}\Z}=\zvf^0(D)X_\Z$,
we have $\lna R,X_\Z\rna_{\wt{\sL}\Z}^{\zvf^0}=X_\Z\we
i_{\zvf^0}R$ for any $R\in\Sec(\sA(\breve{\sL}\Z))$, and, due to
the properties of the Schouten-Jacobi brackets,
\beas
\lna\cJ,\cJ\rna_{\wt{\sL}\Z}^{\zvf^0}&=&
\lna\cJ_0,\cJ_0\rna_{\wt{\sL}\Z}^{\zvf^0}+ 2\lna\cJ_0,X_\Z\we
D_0\rna_{\wt{\sL}\Z}^{\zvf^0}+\lna X_\Z\we D_0 ,X_\Z\we
D_0 \rna_{\wt{\sL}\Z}^{\zvf^0}\\
&=&\lna\cJ_0,\cJ_0\rna_{\wt{\sL}\Z}^{\zvf^0}+ 2(\lna\cJ_0,X_\Z
\rna_{\wt{\sL}\Z}^{\zvf^0}\we D_0-X_\Z\we\lna\cJ_0, D_0
\rna_{\wt{\sL}\Z}^{\zvf^0}-i_{\zvf^0}\cJ_0\we X_\Z\we D_0)\\
&=&\lna\cJ_0,\cJ_0\rna_{{\sL}M}^{\zvf^0}
-2X_\Z\we\lna\cJ_0,X_0\rna_{{\sL}M}^{\zvf^0},
\eeas
that vanishes exactly when
$\lna\cJ_0,\cJ_0\rna_{{\sL}M}^{\zvf^0}=0$ and
$\lna\cJ_0,X_0\rna_{{\sL}M}^{\zvf^0}=0$, i.e. when $\cJ_0$ is a
Jacobi structure for which $D_0$ acts as a derivation of the
corresponding Jacobi bracket
\end{pf}

\medskip
Since aff-Poisson and aff-Jacobi structures have been recognized
as canonical structures, we can apply results of \cite{GU3} to
characterize them in terms of induced morphisms of vector bundles,
and results of \cite{GM1,GM2} to define the corresponding
cohomology and homology.

For $Y\in\Sec(\sA(\wt{\sT}\Z))$ denote by $Y^c$ its Lie algebroid
complete lift to a multivector field on $\wt{\sT}\Z$ (see
\cite{GU1,GU2} or the survey in \cite{GU3}). Similarly, for
$Y\in\Sec(\sA(\wt{\sL}\Z))$ denote by $\wh{Y}_{\zvf^0}$ its Jacobi
algebroid complete lift to a first-order polydifferential operator
on $\wt{\sL}\Z$ (see \cite{GM1} or the survey in \cite{GU3}). Let
$\zL_{\wt{\sT}^\s\Z}$ be the canonical linear Poisson structure on
$\wt{\sT}^\s\Z$ representing the Lie algebroid structure on
$\wt{\sT}\Z$ and let $J_{\wt{\sL}^\s\Z}$ be the canonical
homogeneous Jacobi structure on the dual $\wt{\sL}^\s\Z$ of
$\wt{\sL}\Z$ representing the Jacobi algebroid structure on
$\wt{\sL}\Z$.

\begin{theorem}\
\begin{description}
\item{(i)} $\zL\in\Sec(\bigwedge^2\wt{\sT}\Z)$ represents an aff-Poisson
structure on $\Z$ if and only if the tensors $\zL_{\wt{\sT}^\s\Z}$
and $-\zL^c$ are $\sharp_\zL$-related, where
$\sharp_\zL:\wt{\sT}^\s\Z\ra\wt{\sT}\Z$,
$\sharp_\zL(\zm_m)=i_{\zm_m}\zL(m)$.
\item{(ii)} $\cJ\in\Sec(\bigwedge^2\wt{\sL}\Z)$ represents an aff-Jacobi
structure on $\Z$ if and only if the first-order bidifferential
operators $J_{\wt{\sL}^\s\Z}$ and $-\wh{\cJ}_{\zvf^0}$ are
$\sharp_\cJ$-related, where
$\sharp_\cJ:\wt{\sL}^\s\Z\ra\wt{\sL}\Z$,
$\sharp_\cJ(\zw_m)=i_{\zw_m}\cJ(m)$.
\end{description}
\end{theorem}

\begin{theorem}\
\begin{description}
\item{(a)}  $\zL\in\Sec(\bigwedge^2\wt{\sT}\Z)$ represents an aff-Poisson
structure on $\Z$ if and only if the graded operator
$\pa_\zL(Y)=\lna\zL,Y\rna_{\wt{\sT}\Z}$ of degree 1 on
$\sA(\wt{\sT}\Z)$ is a cohomology operator, i.e. $(\pa_\zL)^2=0$.
\item{(b)} $\cJ\in\Sec(\bigwedge^2\wt{\sL}\Z)$ represents an aff-Jacobi
structure on $\Z$ if and only if the graded operators
$\pa_{\cJ}^t(R)=\lna\cJ,R\rna_{\wt{\sL}\Z}^{\zvf^0}+ti_{\zvf^0}\cJ\we
R$ of degree 1 on $\sA(\wt{\sL}\Z)$ are cohomology operators for
all $t\in\R$, i.e. $(\pa_\cJ^t)^2=0$.
\end{description}
\end{theorem}
\begin{pf} The implication "if" is essentially the graded Jacobi identity
applied to the bracketing with canonical structures. The other
follows from the fact that the corresponding Schouten and
Schouten-Jacobi brackets have no central elements among 3-tensors.
\end{pf}

\medskip\noindent
The cohomology associated to $\pa_\zL$ we will call the {\it
aff-Poisson cohomology}. The cohomology associated to $\pa^0_\cJ$
(resp. $\pa^1_\cJ$) we will call {\it aff-Jacobi cohomology}
(resp., {\it aff-Lichnerowicz-Jacobi cohomology}).

For any Lie algebroid structure on a vector bundle $E$ denote by
$\Ll_Y$ the corresponding Lie differential
$\Ll_Y=i_Y\circ\D_{E}-(-1)^{\vert Y\vert} \D_{E}\circ i_Y$ with
respect to the multisection $Y\in\Sec(\bigwedge ^{\vert
Y\vert}E)$. Here and further $\vert Y\vert$ denotes the degree of
the tensor $Y$.
\begin{theorem} \
\begin{description}
\item{(a)}  $\zL\in\Sec(\bigwedge^2\wt{\sT}\Z)$ represents an aff-Poisson
structure on $\Z$ if and only if the Lie differential
$$\Ll_\zL=i_\zL\circ\D_{\wt{\sT}\Z}-\D_{\wt{\sT}\Z}\circ i_\zL,$$
which is a graded operator of degree -1 on $\sA(\wt{\sT}^\s\Z)$,
is a homology operator, i.e. $(\Ll_\zL)^2=0$.
\item{(b)} $\cJ\in\Sec(\bigwedge^2\wt{\sL}\Z)$ represents an aff-Jacobi
structure on $\Z$ if and only if the Jacobi-Lie differential
$$\Ll^{\zvf^0,t}_\cJ(\zw)=\Ll_\cJ(\zw)+(\vert\zw\vert+t)
i_{i_{\zvf^0}\cJ}(\zw)+\zvf^0\we i_\cJ(\zw),
$$
which is a graded operator of degree -1 on $\sA(\wt{\sL}^\s\Z)$,
is a homology operator for each $t\in\R$, i.e.
$(\Ll^{\zvf^0,t}_\cJ)^2=0$.
\end{description}
\end{theorem}
\begin{pf} The part "if" follows from the identities (see \cite{GM2})
$$2(\Ll_\zL)^2=-\Ll_{\lna\zL,\zL\rna_{\wt{\sT}\Z}}$$
and
$$2(\Ll^{\zvf^0,t}_\cJ)^2=-\Ll^{\zvf^0,t}_{\lna\cJ,\cJ\rna_{\wt{\sL}\Z}^{\zvf^0}}.$$
The other part follows from the fact that passing to the Lie
differentials in the algebroids in question is injective for
3-tensors.
\end{pf}

\medskip\noindent
The homology associated with $\Ll_\zL$ we will call {\it
aff-Poisson homology}. The homology associated with
$\Ll^{\zvf^0,0}_\cJ$ we will call {\it aff-Jacobi homology}.

Aff-Poisson and aff-Jacobi structures give also rise to the
corresponding triangular Lie bialgebroids and Jacobi bialgebroids
(cf. \cite{MX,GM1,IM}). We will not go into the details here.

\section{Applications}
{\bf Example 9.} (\cite{TU,Ur2})  In gauge theories potentials are
interpreted as connections on principal bundles. In the
electrodynamics the gauge group is $(\R, +)$ and the potential is
a connection on a principal bundle $\zz\colon \Z\rightarrow M$
over the space-time $M$, i.e. on an AV-bundle $\Z=(Z, 1_M)$ over
$M$. An electromagnetic potential is a section $\za\colon
M\rightarrow \sP\Z$.

    According to \cite{We}, the phase manifold for a particle with the
charge $e\in \R$ is obtained by the symplectic reduction of
$\sT^{\textstyle *}\Z$ with respect to the coisotropic submanifold
        $$  K_e = \{p\in \sT^{\textstyle *}\Z\colon \langle p, X_\Z\rangle
=-e\}.$$
    Let us denote by $\sP_e\Z$ the reduced phase space. It is easy to see
that it is an affine bundle modelled on $\sT ^{\textstyle *} M$.
We show that $\sP_e\Z$ is the phase bundle for certain special
affine bundle $\bZ_e$.

    First, let $\bY =Z\ati\bI$ be the trivial AV-bundle over $Z$.
    We define an $\R$-action on $\bY$ by the formula
        $$ (Z\times \R)\times \R\ni ((z,r),t) \mapsto (z+t, r +te)\in
Z\times \R =Y.$$
     The space of orbits is an affine bundle modelled on $M\times \R$ and
denoted by $Z_e$. We denote by $\zz_e$ the canonical projection
$Z_e\rightarrow M$. The distinguished section of $\sV(\bY)$ (the
function $ 1_Z$) projects to the constant function $ 1_M$ and the
canonical projection $\zl_e\colon Y\rightarrow Z_e$ is a morphism
of special affine bundles $\bY\rightarrow \bZ_e = (Z_e, 1_M)$. The
induced $\R$-action on $\Z_e$ has the form
        $$ \zl_e(z,r) +s = \zl_e(z,r+s) = \zl_e(z+t, r+ s +te).$$
    For $e=0$ the bundle $\Z_e$ is trivial: $\Z_0 =M\times \R$ and for
$e\neq 0$ we have a diffeomorphism
        $$
    \zF_e\colon Z\rightarrow Z_e, \quad
    \colon z \mapsto \zl_e(z,0).$$
The diffeomorphism $\zF_e$ is an isomorphism of the special affine
bundle $(Z,-\frac{1}{e} 1_M)$ onto $\Z_e$:
        $$ \zF_e(z-\frac{1}{e}r) = \zl_e(z-\frac{1}{e}r,0) = \zl_e(z,r)=
        \zl_e(z,0)+r = \zF_e(z) +r.$$
In particular, $\Z_{-1}=\Z$ and $\Z_1=\ol{\Z}$. To put it simpler,
let us observe that, according to \cite[Example 3]{GGU}, $\Z_e$ is
just the level-set of $1_\Z$ in $\wh{\Z}$ associated with value
$-e$. The diffeomorphism $\zF_e$ comes just from the homotety by
$-e$ in $\wh{\Z}$.

    Let $\zs$ be a section of $\zz_e$. The function $\zl_e^{\textstyle *}
\zs$ on $\Z$ has the property
        $$  X_\Z(\zl_e^{\textstyle *} \zs) =-e.$$

     We conclude that the induced by $\zl_e$ relation $\sP Y\rightarrow
\sP Z_e$ is the symplectic reduction with respect to a coisotropic
submanifold
        $$ K_e = \{p\in \sT^{\textstyle *}\Z\colon \langle p, X_\Z\rangle
=-e\}.$$
    Thus the phase manifold $\sP_e\Z$ for a particle with the charge $e$ is the phase
bundle for the special affine bundle $\bZ_e$. Another way to see
this is to use the decomposition
$\sT^\s\Z=\widetilde{\sT}^\s\Z\ti_M\Z$. The symplectic reduction
in question is the reduction with respect to the moment map for
the phase lift of the canonical $\R$-action on $\Z$, i.e.
$$\sP_e\Z=\{[\za_{z_m}]\in\widetilde{\sT}^\s\Z:
\lan\za_{z_m},X_\Z(z_m)\ran_{z_m}=-e\}.
$$
But $\lan\za_{z_m},X_\Z(z_m)\ran_{z_m}=-e$ is equivalent to
$\lan\za_{z_m},-\frac{1}{e}X_\Z(z_m)\ran_{z_m}=1$ that is a form
of a definition of $\sP\Z_e$. That the symplectic structure on
$\sP\Z$, defined originally as the pull-back from $\sT^\s M$ when
a section of $\Z$ is chosen, coincides with the one reduced from
$\sT^\s\Z$ can be easily checked in the given trivialization.

    The isomorphism $\zF_e$ gives a one-to-one correspondence between
sections of $\zz$ and sections  of $\zz_e$, for $e\neq 0$. It
follows that a chosen section of $\zz$ provides a trivialization
of $\Z$ and also of $\Z_e$. In such trivializations, a section
$\zs$ of $\zz$ and the corresponding section $\zF_e\circ \zs$ of
$\zz_e$ are functions on $M$ related by the formula
        $$ \zF_e\circ \zs (m) = -e\zs(m).$$
     The correspondence $\zs\rightarrow \zF_e\circ \zs$ of sections
projects to a correspondence of affine covectors and consequently
gives a correspondence of affine 1-forms. Let $\za$ be a section
of $\sP \zz:\sP\Z\ra M$ and $\za_e$ be the corresponding  section
of $\sP\zz_e$. In a given trivialization, the sections $\za$ and
$\za_e$ are 1-forms related by the formula $\za_e = -e\za$.

    The lagrangian of a relativistic charged particle is a section $L_e$ of the bundle
$\wt{\sT}\zz_a\colon \wt{\sT}\Z_e\rightarrow \sT M$ over the open
set $C=\{v\in \sT M\colon g(v,v)>0\}$  given by the formula
        $$ L_e(v)= \langle \za_e,v\rangle +  m\sqrt{g(v,v)},$$
     where $g$ is the metric tensor on the space-time $M$, $m$ is the mass of the
particle, and $\langle\za_e,v\rangle=\za_e(v)$, where an element
of $\sP\Z_e$ is interpreted as a linear section of
$\wt{\sT}\zx_e:\wt{\sT}\Z_e\ra\sT M$, i.e. as an element of
$\LS(\wt{\sT}\Z_e)$.
In this example lagrangians are sections of an AV-bundle. Hamiltonians are
ordinary functions but not on a cotangent bundle but on the affine phase
bundle $\sP\Z_e$.

\bigskip\noindent
{\bf Example 10.} (cf. \cite{Ur1}) The space of events for the
inhomogenous formulation of time-dependent mechanics is the space-time $M$
fibrated over the time $\R$. First-jets of this fibration form the
infinitesimal (dynamical) configuration space. Since there is the
distinguished vector field $\pa_t$ on $\R$, the first-jets of the
fibration over time can be identified with those vectors tangent to $M$
which project on $\pa_t$. Such vectors form an affine subbundle $A$ of the
tangent bundle $\sT M$ modelled on the bundle $\sV M$ of vertical vectors.
The Lagrange formalism in the affine formulation originates on the
AV-bundle $A\ati\bI\ra A$ and the lagrangians are ordinary functions on
$A$. The Hamilton formalism now takes place not on the dual vector bundle
$\sV^\s M$ of $\sV M$, as in the classical approach, but on the dual
AV-bundle $\zz:(A\ati\bI)^\#=\bA^\dag\ra\sV^\s M$ which can be recognized
as $\zz:\sT^\s M\ra\sT^\s M/\la dt\ran$ and which carries a canonical
aff-Poisson structure induced from the canonical symplectic Poisson
bracket on $\sT^\s M$ (cf. Example 5). The hamiltonians are sections of
this bundle. To compare with the standard approach, let us assume that we
have a decomposition $M=Q\ti\R$ of the space-time into a product of space
and time. This induces the decomposition $\sT^\s M=\sT^\s Q\ti\sT^\s \R$.
Sections $\zs$ of $\zz$ can be identified with functions (time-dependent
hamiltonians) $H=H(\za,t)$ on $\sV^* M=\sT^\s Q\ti\R$ by
$\zs_H(\za,t)=(\za,t,-H(\za,t)dt)$. The dynamics induced by the section
$\zs_H$ is, as in Example 5, the projection of the dynamics on $\sT^\s M$
induced by $\bF_{\zs_H}$. The distinguished section of $\sT^\s M$ is $dt$,
so that the distinguished section in the AV-bundle $\zz$ is represented by
$-dt$. Thus,
$$\bF_{\zs_H}(\za,t,p)=H(\za,t)+p,$$ where $(t,p)$ are the standard
Darboux coordinates in $\sT^\s \R$. The hamiltonian vector field of
$\bF_{\zs_H}$ on $\sT^\s M$ is therefore $X_{H_t}+\pa_t$, where $X_{H_t}$ is
the hamiltonian vector field of $H_t(x)=H(x,t)$ on $\sT^\s Q$, so we have
recovered the correct dynamics. However, in our picture, the term $\pa_t$ is
not added `by hand' but it is generated from $\zs_H$ by means of the
aff-Poisson structure. Of course, if we have no decomposition into space and
time, there is no canonical $\pa_t$ on $M$ and nothing canonical can be added
by hand in the standard approach. This problem disappears in the aff-Poisson
formulation. In this example, hamiltonians are sections of an AV-bundle and
lagrangians are ordinary functions however not on a vector but on an affine
bundle.

\medskip\noindent
{\bf Example 11.} The last example is devoted to a hamiltonian
formulation of dynamics of one massive particle in the Newtonian
space-time (cf. \cite{GU,Ko}). Even in a fixed inertial frame, up
to now, there was no satisfactory description of the dynamics in
the hamiltonian formulation. First, we would like to present
difficulties that appear while constructing the description for
the dynamics in an inertial frame and then we will show the
solution in the language of AV-geometry.

Let $N$ be the Newtonian space-time i.e. a four-dimensional affine
space equipped with a covector $\zt$ being an element of the dual
of the model vector space $\sV(N)$ and an euclidean metrics $g$ on
the kernel of $\zt$. The covector $\zt$ is used for measuring time
intervals between events and the metrics measures spatial distance
between simultaneous events. We will denote the kernel of $\zt$ by
$E_0$ and the level-1 set of $\zt$ by $E_1$. The vector space
$E_0$ is of course a vector subspace of $\sV(N)$ and $E_1$ is an
affine subspace of $\sV(N)$ modelled on $E_0$. The elements of
$E_1$ are physical velocities of particles. On the other hand,
every element of $E_1$ represents a class of inertial observers
moving in the space-time with the same constant velocity. Such
class of observers will be called an inertial frame. The
configuration space for one massive particle is $N\times E_1$.
Having an inertial frame $u$, we can identify the affine subspace
$E_1$ with its model vector space, the phase space is generally
accepted to be $N\times E_0^\ast$. The correct phase equations for
the potential $\zf$ are:
\bea \label{eqn1}
&\dot x = g^{-1}(\frac{p}{m})+u, \\
\label{eqn2} &\dot p = -\xd_s\zf(x),
\eea
where $(x,p,\dot x, \dot p)$ is an element of $\sT(N\times
E_0^\ast)$ that can be identified with $N\times E_0^\ast\times
\sV(N)\times E_0^\ast$. The subscript in $\xd_s$ means that we
differentiate only in the spatial directions, i.e. vertical with
respect to the projection on time.

The standard hamiltonian description is based on the fact that
$N\times E_0^\ast$ is a Poisson manifold with the Poisson
structure being reduced from the canonical Poisson structure of
$\sT^\ast N\simeq N\times \sV(N)^\ast$. The problem is that from
the hamiltonian
$$
h_u(x,p)=\frac1{2m}\langle p, g^{-1}(p)\rangle +\zf(x)
$$
we obtain the vector field which is vertical with respect to the
projection on time:
\beas
& \dot x = g^{-1}(\frac{p}{m}), \\
& \dot p = -\xd_s\zf(x). \eeas
    Any vertical vector field cannot be a physical
motion, so we have to add `by hand' the constant vector field $u$. As in the
previous example, this problem can be solved by replacing the Poisson
structure on the phase manifold by an affine Poisson structure. However, the
equations of motion as well as the affine Poisson structure depend on the
choice of the reference frame.

    To get frame-independent formulation for the dynamics, let us consider
first frame-dependent lagrangian $\ell_u$. It is a function defined on
$N\times E_1$
$$
\ell_u(x,v)=\frac{m}{2}\langle g(v-u), v-u \rangle-\varphi(x).
$$
Let us look at the solution of this problem. If $u$ and $u'$ are
two inertial frames then the difference
$$f_{u,u'}(v)=\ell_u(x,v)-\ell_{u'}(x,v)=m\langle g(u'-u),
v-\frac{u+u'}{2}\rangle
$$
is an affine function on $E_1$. Now we define an equivalence
relation $\sim_{\ell}$ in the set $E_1\times E_1\times \R$ by
$$
(u,v,r)\sim_{\ell}(u',v',r')\quad\Longleftrightarrow\quad
v=v',\,\, r=r'+f_{u,u'}(v).
$$
The set of equivalence classes for $\sim_{\ell}$ will be denoted by $A_0$.
We observe that since $f_{u,u'}$ is an affine function, $A_0$ is an affine
space of dimension $4$. There is a projection from $A_0$ to $E_1$. The
model vector space for $A_0$ is $(E_1\times E_0\times \R)\slash
\sim_{v\ell}$, where the equivalence relation $\sim_{v\ell}$ is in a sense
the linear part of $\sim_\ell$: we say that two elements $(u,w,r)$ and
$(u',w',r')$ of $E_1\times E_0\times \R$ are equivalent if
$$w=w',\,\,\, r=r'+m\langle g(u-u'),w\rangle. $$
In $\sV(A_0)$ we distinguish an element $w_0=[u,0,1]$ so now
$\bA_0=(A_0,w_0)$ is a special affine space and $\bA=N\times\bA_0$
is a special affine bundle over $N$. The mapping
$$
(x,v)\longmapsto (x,[u,v,\ell_u(x,v)])
$$
is a section of $\bAP(\bA)$ and can be understood as
frame-independent lagrangian. Note that it is no longer a function
but a section of an AV-bundle. Any section of $\bA$ can be
represented in the form
$$ x\longmapsto \tilde X(x)=[u,X(x), r(x)],$$
where $X$ is a vector field on $N$ with values in $E_1$ and $r$ is
a function on $N$. We define a bracket on sections of $\bA$ by the
following formula
$$
[\tilde X, \tilde Y]=[u,[X,Y], Xs-Yr],
$$
where $\tilde X(x)=[u,X(x), r(x)]$ and $\tilde Y(x)=[u,Y(x),
s(x)]$. The definition is correct. Indeed, if we have other
representatives $X(x)=[u',X(x),r(x)-f_{u,u'}(X(x))]$ and
$Y(x)=[u',Y(x),s(x)-f_{u,u'}(Y(x))]$, then, since $f_{u,u'}$ is
affine, $(\Ll_{X(x)}(f_{u,u'}\circ
Y))(x)=(f_{u,u})_\sv\circ(\Ll_{X(x)}Y)(x)$, where $\Ll_{X(x)}$ is
the directional derivative in the direction $X(x)$ and
$(f_{u,u'})_\sv$ is the vector part of $f_{u,u'}$. Moreover,
$(\Ll_{X(x)}Y-\Ll_{Y(x)}X)(x)=[X,Y](x)$, so that we get
$$X(r-f_{u,u'}\circ Y)-Y(s-f_{u,u'}\circ X)=X(r)-Y(s)-
(f_{u,u'})_\sv\circ[X,Y],$$ that proves the correctness of the
definition.

Having two vector fields $X,Y$ with values in $E_1$ we have
$$0=(\xd \tau)(X,Y)=X\langle\tau, Y\rangle-Y\langle \tau, X\rangle
+\langle\tau,[X,Y]\rangle.$$ Since $X\langle\tau,
Y\rangle=Y\langle \tau, X\rangle=0$ we obtain that
$\langle\tau,[X,Y]\rangle=0$, i.e. $[X,Y]\in E_0$. The bracket of
two sections of $\bA$ is therefore a section of $\sV(\bA)$ and it
is easy to see that it is a Lie affgebroid bracket with the anchor
morphism $\gamma: \A\longrightarrow \sT N$ defined as
$\gamma([u,X,r])=X$. Moreover, the section $w_0$ is central for
the bracket, i.e. $[\wt{X},w_0]^2_\sv=0$ for all $X$. Therefore,
according to the Theorem \ref{c3}, we have that the corresponding
aff-Jacobi bracket on the AV-bundle ${\bAP(\bA^\#)}$ is
aff-Poisson. We claim that this structure is the correct structure
for generating the equations of motion for the hamiltonian
formulation of the dynamics in question.

Indeed, $\bA^\#$ is by definition $\bAff(\bA,\bI)$. Like in the
lagrangian case, we will represent it as the set of cosets of an
appropriate equivalence relation. In the space $E_1\times
E_0^\ast\times\R$ we define an equivalence relation $\sim_h$ by
$$(u,p,s)\sim_h(u',p',s')\quad\Longleftrightarrow\quad
p=p'+mg(u-u'),\,\,\, s=s'+\langle p,u-u'\rangle+\frac12 m\langle
g(u-u'), u-u' \rangle.
$$
An equivalence class $[u,p,s]$ represents an affine function
$\zx_{[u,p,s]}$ on $A_0$ given by
$$
\zx_{[u,p,s]}([u,v,r])=\langle p, v-u\rangle + s - r.
$$
Its linear part $(\zx_{[u,p,s]})_{\sV}([u,w,r])=\langle p,w\rangle
- r$ gives $-1$ while evaluated on $w_0$. The model vector space
for $\bA_0^\#$ is $\sV(N)^\ast$ with distinguished element $\zt$.
We have:
$$ [u,p,s]+\zp=[u, p+\imath(\zp), s+\langle \zp, u\rangle],$$
where $\imath$ is the canonical projection from $\sV(N)^\ast$ on
$E_0^\ast$. Let us denote by $P$ the space of affine momenta, i.e.
the space $E_0\times E_0^\ast\slash \sim_P$ for the relation
$$(u,p)\sim_P(u',p')\quad\Longleftrightarrow\quad p=p'+mg(u-u').$$
We observe that $\bAP(\bA^\#)$ is an AV-bundle over $N\times P$.
The frame-independent hamiltonian is a section of $\bAP(\bA^\#)$:
$$
(x, [u,p])\longmapsto h((x, [u,p])=(x, [u,p, h_u(x,p)]).
$$
Using the canonical aff-Poisson structure on $\bAP(\bA^\#)$ we can
generate out of $h$ an affine derivation of $\bAP(\bA^\#)$, i.e. the
section of $\wt{\sT}(\bAP(\bA^\#))$. This section projects to a vector
field on $N\times P$ that is understood as the equation of motion.

Now, let us calculate the equations of motion in coordinates. For,
we choose an inertial frame $u$ and the coordinates $(x^0, x^i)$,
$i=1,2,3$ such that $\partial_0=u$. By $(p_i)$ we denote the
adapted coordinates on $E_0^\ast$. Using the inertial frame we
have the following identifications:
\beas
\bA_0 & \simeq & E_1\times \bI, \\
\bA & \simeq & N\times E_1\times \bI, \\
\bA^\# & \simeq & N\ti E_0^\ast\times \bI, \\
\sss{Sec}(\bAP(\bA^\#)) & \simeq & {\cal C}^\infty(N\times E_0^\ast).
\eeas
The bracket of sections of $N\times E_1\times \R,$ has the obvious
form $[(X, \zx),(Y, \zy)]=([X,Y], X\zy-Y\zx)$. If the sections
take values in $E_1$ then the bracket takes values in $E_0$. The
affine function on $N\times E_0^\ast$ that corresponds to the
section $(X,\zx)$ is
$$\imath_{(X,\zx)}(x,p)=\langle p, X \rangle-\zx.$$
The Poisson bracket for functions corresponding to sections $(X,
\zx)$, $(Y, \zy)$ is given by the formula
$$
\{\imath_{(X,\zx)},\imath_{(Y,\zy)}\}=\imath_{[(X,\zx),(Y,\zy)]}=
\langle p, [X,Y] \rangle-X\zy+Y\zx,
$$
which in coordinates reads
$$
\{\imath_{(X,\zx)},\imath_{(Y,\zy)}\}=
p_iX^j\partial_jY^i-p_iY^j\partial_jX^i+
p_i\partial_0Y^i-p_i\partial_0X^i-X^i\partial_i\zy-
\partial_0\zy+Y^i\partial_i\zx+\partial_0\zx.
$$
From the above formula we obtain that
$$
\{h,\cdot\}=(\partial^j h)\partial_j+\partial_0-(\partial_j
h)\partial^j -\partial_0h.
$$
The vector part of the above operator is exactly what we had in
({\ref{eqn1}}) and ({\ref{eqn2}}). In this example both, hamiltonians and
lagrangians are sections of AV-bundles and not ordinary functions.

\bigskip\noindent
%%%%%%%%%%%%%%%%%%%%%%%%%%%%%%%%%%%%%%%%%%%%%%%%%%%%%%%%%%%%%o

\end{document}